\documentclass[a4paper,fleqn]{cas-dc}

\usepackage[numbers,sort&compress]{natbib}

\usepackage{bm}
\usepackage{comment}
\def\tsc#1{\csdef{#1}{\textsc{\lowercase{#1}}\xspace}}
\tsc{WGM}
\tsc{QE}
\tsc{EP}
\tsc{PMS}
\tsc{BEC}
\tsc{DE}

\usepackage{framed} 
\usepackage{multicol} 
\usepackage{nomencl} 
\makenomenclature
\renewcommand*\nompreamble{\begin{multicols}{2}}
\renewcommand*\nompostamble{\end{multicols}}

\usepackage{algorithm}
\usepackage{algpseudocode}
\usepackage{amsmath,amsfonts,amsthm}
\usepackage{graphicx}
\usepackage{subcaption}
\usepackage{url}
\usepackage{hyperref} 
\usepackage{xcolor}
\usepackage{makecell}

\hypersetup{
    colorlinks=true, 
    linkcolor=black, 
    filecolor=black, 
    urlcolor=black, 
    citecolor=black, 
    linkbordercolor={1 1 1}, 
    pdfborder={0 0 0} 
}

\begin{document}
\let\WriteBookmarks\relax
\def\floatpagepagefraction{1}
\def\textpagefraction{.001}
\let\printorcid\relax 
\shorttitle{\rmfamily Real-time Hosting Capacity Assessment for Electric Vehicles}
\shortauthors{\rmfamily Yingrui~Zhuang et~al.}

\title [mode = title]{Real-time Hosting Capacity Assessment for Electric Vehicles: A Sequential Forecast-then-Optimize Method}

\affiliation[1]{organization={State Key Laboratory of Power System and Generation Equipment, Department of Electrical Engineering, Tsinghua University},
                city={Beijing},
                postcode={100084}, 
                country={China}}

\affiliation[2]{organization={Department of Earth and Environmental Engineering, Columbia University},
                city={New York},
                postcode={10027}, 
                country={United States}}

\affiliation[3]{organization={Department of Mechanical and Automation Engineering, The Chinese University of Hong Kong},
                city={Hong Kong},
                country={China}}

\author[1]{Yingrui~Zhuang}[style=chinese]
\cormark[1]
\ead{zyr21@mails.tsinghua.edu.cn} 
\credit{Conceptualization, Methodology, Software, Writing}

\author[1]{Lin~Cheng}[style=chinese]
\credit{Conceptualization, Supervision}

\author[2]{Ning~Qi}[style=chinese]
\credit{Methodology, Data curation, Reviewing and Editing}

\author[1]{Xinyi~Wang}[style=chinese]
\credit{Data Curation, Software, and Writing}

\author[3]{Yue~Chen}[style=chinese]
\credit{Conceptualization, Methodology, Writing – Reviewing and Editing}

\cortext[cor1]{Corresponding author}

\begin{abstract}
Hosting capacity (HC) assessment for electric vehicles (EVs) is crucial for EV secure integration and reliable power system operation.
Existing methods primarily focus on a long-term perspective (e.g., system planning), 
and consider the EV charging demands as scalar values, 
which introduces inaccuracies in real-time operations due to the inherently stochastic nature of EVs.
In this regard, this paper proposes a real-time HC assessment method for EVs through 
a three-step process, involving real-time probabilistic forecasting, risk analysis and probabilistic optimization.
Specifically, we conduct real-time probabilistic forecasting to capture the stochastic nature of EV charging demands across multiple charging stations by performing deterministic forecasting and fitting the distribution of forecasting errors.
The deterministic forecasting is conducted
using an adaptive spatio-temporal graph convolutional network (ASTGCN).
ASTGCN leverages adaptive spatial feature extraction, attention-based temporal feature extraction, 
and second-order graph representation to improve the forecasting performance.
Subsequently, based on the probabilistic forecasting of EV charging demands,
we conduct real-time risk analysis and operational boundary identification by
utilizing probabilistic power flow calculations to assess potential violations of secure operation constraints.
Furthermore, we present the formulation of real-time HC of EVs
considering expected satisfaction of stochastic EV charging demands,
and propose an optimization model for real-time HC assessment of EVs.
Numerical experiments on a real-world dataset demonstrate that the proposed ASTGCN model outperforms state-of-the-art forecasting models
by achieving the lowest root mean square error of 0.0442,
and the real-time HC is improved by 64\% compared to long-term HC assessment.
\end{abstract}

\begin{highlights}
    \item Spatio-temporal graph convolutional network for multiple charging demands forecasting.
    \item Adaptive time-invariant and time-varying spatial feature extraction.
    \item Real-time risk analysis via Gaussian mixture-based probabilistic power flow.
    \item Optimization model for real-time hosting capacity assessment for electric vehicles.
\end{highlights}


\maketitle

\renewcommand{\nomgroup}[1]{%
	\item[
      \textbf{%
        \ifthenelse{\equal{#1}{A}}{\textit{Abbreviations}}{}%
        \ifthenelse{\equal{#1}{B}}{\textit{Sets}}{}%
        \ifthenelse{\equal{#1}{C}}{\textit{Parameters and Variables}}{}%
        \ifthenelse{\equal{#1}{E}}{\textit{Functions}}{}%
	}]%
}
\nomenclature[A]{EV}{Electric vehicle}
\nomenclature[A]{HC}{Hosting capacity}
\nomenclature[A]{ASTGCN}{Adaptive spatio-temporal graph convolutional network}
\nomenclature[A]{GMM}{Gaussian mixture model}
\nomenclature[A]{MAE}{Mean absolute error}
\nomenclature[A]{RMSE}{Root mean square error}
\nomenclature[A]{WAPE}{Weighted absolute percentage error}

\nomenclature[B]{$\varOmega^\mathrm{L}$}{Set of load buses}
\nomenclature[B]{$\varOmega^\mathrm{C}$}{Set of charging stations buses}
\nomenclature[B]{$\mathcal{G}/\mathcal{V}/\mathcal{E}$}{Graph/node/edges}
\nomenclature[B]{$\mathcal{D}$}{Training dataset}

\nomenclature[C]{$N^\mathrm{C}$}{Number of charging stations}
\nomenclature[C]{$N^\mathrm{F}$}{Number of deterministic forecasting intervals}
\nomenclature[C]{$T$}{Historical time slots for forecasting}
\nomenclature[C]{$\bm{X}$}{Feature matrix}
\nomenclature[C]{$\tilde{\bm{L}}$}{Normalized graph Laplacian matrix}
\nomenclature[C]{$\bm{W}$}{Adaptive weighted adjacency matrix}
\nomenclature[C]{$\bm{W}^\mathrm{TI}/\bm{W}^\mathrm{TV}$}{Adaptive time-invariant/time-varying weighted adjacency matrix}
\nomenclature[C]{$\bm{E}$}{Embedding matrix of $\bm{W}^\mathrm{TI}$}
\nomenclature[C]{$\bm{M}$}{Learnable matrix for Mahalanobis distance}
\nomenclature[C]{$\hat{P}_{i,t}/\hat{\mathcal{P}}_{i,t}$}{Deterministic/probabilistic forecasting of EV charging demands at station $i$ and time slot $t$}
\nomenclature[C]{$P_{i,t}$}{True value of EV charging demands at station $i$ and time slot $t$}
\nomenclature[C]{$\bm{A^{\mathrm{tt}}}$}{Temporal attention matrix}
\nomenclature[C]{$\bm{J}/\bm{S}$}{Jacobian/sensitivity matrix of power flow equations}
\nomenclature[C]{$\pi_{i,k}^{(j)}/\mu_{i,k}^{(j)}/\sigma_{i,k}^{(j)}$}{Weight/mean/standard deviation of the $k$-th Gaussian component for station $i$ in the $j$-th forecasting interval}
\nomenclature[C]{$\bm{Z}$}{Bilinear mapping matrix in second-order pooling}
\nomenclature[C]{$\theta_k$}{Chebyshev polynomial coefficient of order $k$}
\nomenclature[C]{$\bm{Q}/\bm{K}/\bm{V}$}{Query/key/value matrix of temporal attention mechanism}
\nomenclature[C]{$V_i/\overline{V}_i/\underline{V}_i$}{Voltage, upper/lowervoltage limit at bus $i$}
\nomenclature[C]{$r_{ij}/x_{ij}$}{Resistance/reactance of line $ij$}
\nomenclature[C]{$P_{ij}^{\mathrm{Ln}}/Q_{ij}^{\mathrm{Ln}}$}{Active/reactive power flow of line $ij$}
\nomenclature[C]{$I_{ij}/\overline{I}_{ij}$}{Current/upper current magnitude of line $ij$}
\nomenclature[C]{$p_i/q_i$}{Active/reactive power injection at bus $i$}
\nomenclature[C]{$P_i^L/Q_i^L$}{Active/reactive load demand at bus $i$}
\nomenclature[C]{$\overline{P}_{i,t}$}{Maximum accommodated EV charging demands at station $i$ and time slot $t$}

\nomenclature[E]{$T_k(\cdot)$}{Chebyshev polynomial of order $k$}
\nomenclature[E]{$\psi(\cdot)$}{Activation function}
\nomenclature[E]{$d^\mathrm{M}(\cdot)$}{Generalized Mahalanobis distance function}
\nomenclature[E]{$\odot$}{Element-wise product}
\nomenclature[E]{$\mathcal{P}/\mathcal{F}$}{Probability density function/cumulative distribution function}
\nomenclature[E]{$\mathrm{MLP}(\cdot)$}{Multi-layer perceptron}

\begin{table*}[!t]
\begin{framed}
    {
        \fontfamily{ptm}\selectfont 
        \printnomenclature
    }
\end{framed}
\end{table*}

\section{Introduction}
With the rapid advancement of electrification in the transportation sector, 
electric vehicles (EVs) have quickly proliferated and formed a deep integration with power systems~\cite{EVimpact}.
According to the report from the International Energy Agency~\cite{iea2024EV}, 
in 2024, EV sales could reach around 17 million.
Besides, in China, EVs' market share could reach up to 45\% of total car sales, 
with their charging demands accounting for more than 30\% of urban residential electricity load demand.
However, the massive integration of EVs has significantly increased the electricity demand,
posing potential risks to power system operations~\cite{impact1,OD}, 
such as overvoltage, network thermal overloading, increasing load peak-valley differences, network power loss, 
and power quality issues.
Along with large-scale transitions, 
existing power system infrastructure is challenged to accommodate the massive integration of EVs
as safely as possible towards a zero-carbon system~\cite{xie2021toward} .
Researchers have defined the concept of hosting capacity (HC)~\cite{HC,HC_review,HC2,HC3} 
to quantify the maximum EV charging demands that can be accommodated by a power system without violating operational constraints 
and incurring substantial upgrades~\cite{zero_corbon}.
To ensure the secure integration of EVs and achieve the goal of net zero carbon emissions, 
it is critical to assess the HC of EVs.

Existing HC assessment methods for EVs primarily focus on a long-term perspective (e.g., system planning),
and consider the EV charging demands as scalar values~\cite{HCEVreview}. 
However, in real-time operation, 
EV charging demands show substantial uncertainty due to stochasticity in user behavior 
and other external factors (e.g., weather conditions, electricity charging price)~\cite{EVbehaviour}, 
and introduce considerable stochastic risks~\cite{EVimpact1} to power system operation.
In this case,
long-time HC estimation with scalar modeling of EV charging demands 
cannot reliably capture the stochastic uncertainty, 
potentially leading to inaccuracies in risk analysis and HC assessment,
and cannot be reliable for decision-making in real-time operation.
Recognizing this, 
Ref.~\cite{DHC} highlights the importance of conducting HC assessments across different time periods, 
utilizing time series data (e.g., forecasted demand values) to more accurately represent uncertain conditions.
Therefore, 
it is essential to develop real-time HC assessment to ensure 
the secure integration of EVs and the reliable operation of power systems.

\begin{figure}[tb]
    \setlength{\abovecaptionskip}{-0.1cm}
    \setlength{\belowcaptionskip}{-0.1cm}
    \centerline{\includegraphics[width=0.95\columnwidth]{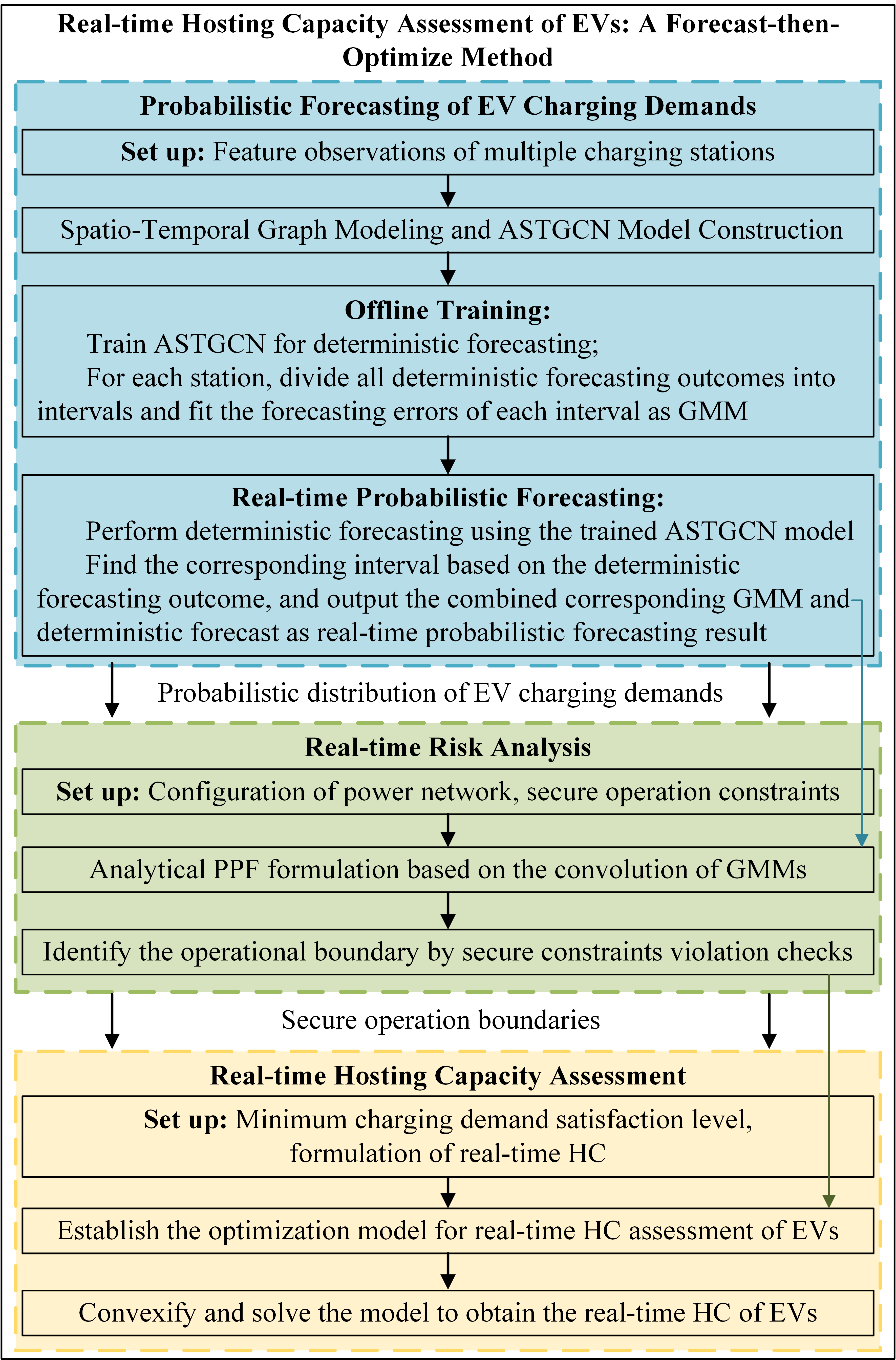}}
    \caption{\rmfamily Flowchart of the proposed real-time HC assessment method for EVs.}
    \label{fig:flowchart}
    \vspace{-0.5cm}
\end{figure}

In this paper, we propose a forecast-then-optimize method for real-time HC assessment for EVs. 
A flowchart is provided in Fig.~\ref{fig:flowchart}.
Specifically, the main contributions are threefold:

\begin{enumerate}
    \item \textbf{Probabilistic forecasting:}
        We propose a probabilistic forecasting method to characterize the stochastic nature of EV charging demands,
        which is achieved through deterministic forecasting and distribution fitting of forecasting errors based on predefined forecasting intervals.
        Specifically, an adaptive spatio-temporal graph convolutional network (ASTGCN) is proposed for 
        deterministic forecasting of EV charging demands across multiple charging stations.
        Compared to existing methods, 
        ASTGCN leverages an adaptive graph convolutional layer to extract the underlying time-invariant and time-varying spatial features of the charging demands.
        Subsequently, we incorporate a spatio-temporal graph convolutional block that 
        combines adaptive spatial feature extraction with attention-based temporal feature extraction to comprehensively capture the spatio-temporal characteristics,
        and further introduces a second-order graph representation to enhance the model's graph representation capability.
        Simulation results indicate that ASTGCN achieves superior performance across multiple evaluation metrics compared to state-of-the-art forecasting models.
    \item \textbf{Risk analysis:}
        We propose a real-time risk analysis and operational boundary identification method.
        Specifically, 
        with probabilistic forecasting the stochastic EV charging demands as Gaussian mixture model (GMM),
        we develop an analytical formulation of probabilistic power flow (PPF) through convolution of GMMs. 
        Then, real-time risk is performed,
        and secure operational boundaries are identified by assessing the potential violations of secure operation constraints.
    \item \textbf{HC assessment:} 
        Considering the uncertain characteristics of EV charging demands in real-time operation,
        we present the formulation of real-time HC of EVs as the maximum expected EV charging demands that can be accommodated by the power system without violating operational constraints.       
        Subsequently, we propose an optimization model to assess the real-time HC of EVs.
        Compared to existing long-term HC assessment methods, the proposed real-time HC assessment method 
        improves the expected EV charging demands accommodation by 64\%, making it better suited to stochastic real-time operation.
\end{enumerate}

The remainder of the paper is organized as follows. 
Section~\ref{Literature Review} reviews the previous works on EV charging demands forecasting, risk analysis and HC assessment.
Section~\ref{ProbASTGCN} introduces the detailed procedure of probabilistic forecasting and the structure of ASTGCN model.
Section~\ref{Probabilistic Risk Analysis} presents the real-time risk analysis method utilizing GMM-based PPF.
Section~\ref{Real-time Hosting Capacity Assessment for EVs} proposes the real-time HC assessment method for EVs.
Numerical studies based on real-world data are provided in Section~\ref{case study} to illustrate comparative performance. 
Finally, conclusions are summarized in Section~\ref{conclusion}.

\section{Literature Review}\label{Literature Review}
\subsection{EV Charging Demands Forecasting}
In terms of forecasting formulation, existing literature can be categorized as deterministic forecasting and probabilistic forecasting.
Deterministic forecasting generates expected values of EV charging demands~\cite{EV_tcn},
while probabilistic forecasting provides statistical information in the form of density~\cite{EVtf}, 
quantiles~\cite{prob_pred_EV2}, or intervals~\cite{interval}.
Probabilistic forecasting methods can be further categorized into parametric methods and non-parametric methods.
Parametric methods assume a specific probability distribution (e.g., normal distribution, beta distribution) for the forecast target,
and directly estimate the parameters of the assumed probability distribution~\cite{prob_pred_PV}.
Parametric methods are suitable for data with a clear distribution. 
Non-parametric methods (e.g., quantile regression~\cite{prob_pred_EV2}, kernel density estimation~\cite{kernal}), 
do not assume a predefined distributional assumptions. 
Instead, they directly model the distribution based on the data itself~\cite{non_para_pred}, 
offering more flexibility to capture complex, irregular patterns,  
but at the cost of increased computational complexity.
Additionally, probabilistic forecasting can be conducted by directly generating the probabilistic distribution of the forecasting target,
or by superimposing a probabilistic error distribution onto deterministic forecasts~\cite{Versatile}. 
Existing literature on EV charging demands forecasting mainly focuses on deterministic forecasting.
However, 
due to the inherent uncertainty in real-time EV charging behavior, 
deterministic forecasting falls short for
overlooking this stochastic nature~\cite{point_pred_not_reliable}. 
Probabilistic forecasting addresses this gap by quantifying forecast uncertainty, 
thereby enabling more robust decision-making in power systems facing uncertainty.

In terms of specific forecasting methods, existing research can be categorized into statistical methods and machine learning methods.
In statistical methods, EV charging behaviors (e.g., charging start time, charging duration, average charging power) 
are mainly simulated using self-defined statistical distributions based on various influencing factors~\cite{EVbehaviour,tripchain},
such as physical parameters of EVs~\cite{EV_characteristic} (e.g., EV type, battery parameters, state of charge),
traffic flow~\cite{EVtf} (e.g., road congestion, driving distance), 
external environmental conditions~\cite{weather} (e.g., temperature), 
and characteristics of charging facilities (e.g., charging station location, available piles, charging price). 
The relationship between these factors and EV charging behaviors can be derived from historical observations or travel surveys~\cite{survey}.
Then, Monte Carlo simulation~\cite{MC1}, origin-destination analysis~\cite{OD}, queuing theory~\cite{EVquene}, and Markov chain~\cite{markovchain} 
are used in previous work to generate EV charging scenarios. 
Additionally, Historical Average (HA) and Autoregressive Integrated Moving Average (ARIMA)~\cite{ARIMA} are commonly employed
for forecasting using temporal changing patterns.
However, statistical methods have general limitations in three aspects:
\textit{(i)} 
incorporating diverse real-world factors into statistical models is challenging, 
and their impacts on EV behavior vary widely across different situations~\cite{survey},
thus facing limitations in scalability and generalization in practice;
\textit{(ii)} substantial diverse data is required for model validation, 
which is often difficult to obtain due to privacy concerns;
\textit{(iii)} statistical regression methods are not sufficient to extract complex coupling features between charging stations.

In recent years, advancements in big data analytics and the growing data accumulation 
have made machine learning methods new opportunities to analyze EV charging behaviors~\cite{PCMP2}.
Machine learning models facilitate end-to-end forecasting by uncovering complex, hidden patterns within extensive datasets.
Existing works have investigated the performance of k-nearest neighbor~\cite{EV_pred1}, 
temporal convolutional network (TCN)~\cite{EV_tcn}, 
convolutional neural network (CNN)~\cite{EVtf},
long short-term memory (LSTM)~\cite{EV_lstm} and gated recurrent units~\cite{EV_pred3}.
These methods generally frame EV charging demands forecasting as a time-series problem, focusing primarily on temporal correlations.
Ref~\cite{STMGCN} advances this by capturing both spatial and temporal correlations, 
proposing a multi-graph spatio-temporal convolutional network (STGCN) for deterministic forecasting of EV charging demands. 
This model considers the Euclidean distance between charging station locations and historical demands to represent  spatial relationships. 
However, given that EV charging demands are coupled with both the power network and transportation network, 
Euclidean distance alone is insufficient to capture the complex spatial correlations.

Overall speaking,
existing methods for EV charging demands forecasting primarily focus on deterministic forecasting,
limiting their ability to account for the inherent stochasticity and uncertainty in real-time EV charging behaviors. 
Furthermore, 
the prevalent use of single-model approaches often overlooks critical spatial coupling characteristics, 
constraining their capacity to effectively capture the complex spatio-temporal dynamics that underpin interactions within power and transportation networks.

\subsection{Risk Analysis}

The significant randomness and uncertainty of real-time EV charging demands
introduce considerable operational risks to power systems~\cite{EVimpact1}.
Generally, risk is quantified 
based on a probabilistic description of the power system state
obtained through probabilistic power flow calculations~\cite{PPFreview}.
Simulation methods (e.g., Monte Carlo simulation), analytical methods (e.g., cumulant transformation), and approximation methods (e.g., point estimation)
are commonly used for PPF calculation.
In simulation methods,
numerous discrete EV integration scenarios are generated using Monte Carlo methods~\cite{PVHC1,PVHC2},
then power flow is calculated for each scenario to calculate the probabilistic distribution of system state.
However, numerous power flow calculations make simulation methods computationally demanding.
Analytical methods and approximation methods assume that the PPF follows a quasi-Gaussian distribution, 
and reconstruct the PPF from low-order statistical quantities of random variables using series expansion or maximum entropy methods. 
Despite their efficiency, these methods are limited by the Gaussian assumptions, 
which often do not hold in real-time operations where power flow distributions can deviate significantly from Gaussian forms, 
even often resulting in multimodal distributions due to complex, dynamic system behaviors (e.g., load shedding strategies altering the distribution of EV charging demands).
Ref.~\cite{PPFGMM} employs GMM into load flow formulation with a simplified backwards-forwards load flow calculation.

\subsection{HC Assessment}
Existing HC assessment methods primarily focus on a long-term perspective (e.g., system planning)
and model the EV charging demands as scalar values~\cite{HCEVreview}.
Existing methods can be generally classified into simulation methods and optimization methods.

Simulation methods generate numerous scenarios using Monte Carlo simulation~\cite{HCMC}
to describe the locations of EV charging stations and the EV charging demand level.
For each scenario, power flow calculations are performed to compute network-related variables, 
such as voltages and power flows, and to assess potential violations of operational constraints.
In some studies, sensitivity analysis is applied by gradually increasing the EV penetration level until system secure operational constraints
are violated.
After calculating numerous scenarios, HC is determined as the highest feasible level of EV integration
that still satisfies all system secure operational constraints.
While certain studies have proposed accelerated power flow calculation techniques to 
improve computational efficiency~\cite{pf_fast1,pf_fast4},
scenario-based simulation methods remain computationally intensive, 
especially in large networks, making them impractical for real-time applications.

Optimization-based HC assessment methods formulate HC assessment as a maximization optimization problem, 
where EV locations and power injection levels are treated as decision variables 
to maximize EV integration within system secure operational limits~\cite{EVHC1,HCoptimization,HC_tonglang}.
Some studies propose to further improve HC 
by incorporating the collaborative operation of EVs with other distributed energy resources 
to reduce overall system operational risks~\cite{EVoperation3,EVoperation4}.
However, these methods consider EV charging demands as scalar values,
neglecting the inherent stochastic variability in real-time EV charging behavior, 
which can lead to inaccuracies and reliability issues.
To ensure secure EV integration and reliable power system operation, 
it is crucial to develop real-time HC assessment methods that incorporate the stochastic nature of EV charging demands.

\section{Probabilistic Forecasting of Multiple EV Charging Demands}\label{ProbASTGCN}

In this section, 
we first introduce the overall procedures of the proposed probabilistic forecasting method for multiple EV charging demands,
and then present the detailed structure of the proposed adaptive spatio-temporal graph convolutional network.

\subsection{Procedures of Real-time Probabilistic Forecasting}

The real-time probabilistic forecasting task for multiple EV charging demands can be defined as follows: 
at time slot $t$, based on a historical feature matrix $\bm{X}_t$ from the previous $T$ time slots for $N^{\mathrm{C}}$ charging stations in charging stations set $\varOmega^{\mathrm{C}}$, 
represented as $\bm{X}_t = \{ X_{i, \tau} \mid i \in \varOmega^{\mathrm{C}}, \tau \in [t - T, t - 1] \} \in \mathbb{R}^{N^{\mathrm{C}} \times T}$,
we aim to forecast the probabilistic distributions of EV charging demands of the next time slot,
which is denoted as $\hat{\bm{\mathcal{P}}}_{t} = \{ \hat{\mathcal{P}}_{i, t} \mid i \in \varOmega^{\mathrm{C}}
\}$. 
We have $|\varOmega^{\mathrm{C}}|=N^{\mathrm{C}}$.
For clarity in the following sections, 
the notations $ \mathcal{P} $, $ P$, $ \bm{P}$, used throughout this paper will refer specifically to EV-related variables, 
unless otherwise indicated by special superscripts or annotations.

In this paper, 
we achieve probabilistic forecasting of multiple EV charging demands
through two stages: offline training and real-time forecasting.
The detailed procedures are as follows:

\subsubsection{Stage 1: Offline Training}

\textit{\textbf{Step 1.1:}} 
A deterministic forecasting model is first developed and trained.
Specifically, we utilize an adaptive spatio-temporal graph convolutional network to achieve deterministic forecasting.

\textit{\textbf{Step 1.2:}} 
Deterministic forecasting is conducted for each charging station over the past $T+1$ to $t-1$ time slots.
$\hat{P}_{i,\tau}$ denotes the deterministic forecasting of the $i$-th charging station at time slot $\tau$,
with $i \in \varOmega^{\mathrm{C}}, \tau \in [T+1, t-1]$.
The corresponding true value is denoted as $P_{i,\tau}$.

\textit{\textbf{Step 1.3:}}  
The forecasting zone (i.e., $\left[0,1\right]$) is divided into $N^{\mathrm{F}}$ intervals with length $1/N^{\mathrm{F}}$.
The $j$-th interval is denoted as 
$I^{\mathrm{F}}_j=\left[(j-1)/N^{\mathrm{F}},j/N^{\mathrm{F}}\right]$.
Note that more forecast intervals can be built with sufficient data accumulation.
Smaller interval length leads to more accurate probabilistic forecasting results.
For each station $i$,
we collect all the deterministic forecasting outcomes at all previous time slots that fall into the $j$-th interval,
i.e., $\hat{\bm{P}}_{i}^{(j)}=\{ \hat{P}_{i,\tau}|\hat{P}_{i,\tau}\in I^{\mathrm{F}}_j,\ \tau \in [T+1,t-1]\}$,
and the corresponding forecasting errors are denoted as $\bm{\tilde{P}}_{i}^{(j)}=\{P_{i,\tau}-\hat{P}_{i,\tau}|\hat{P}_{i,\tau}\in I^{\mathrm{F}}_j\}$.

\textit{\textbf{Step 1.4:}}  
$\bm{\tilde{P}}_{i}^{(j)}$ is fitted into probabilistic distributions utilizing Gaussian mixture model:
\begin{equation}\label{GMMforecast}
    \bm{\tilde{P}}_{i}^{(j)} \sim \tilde{\mathcal{P}}_{i}^{(j)}=\sum_{k=1}^{K_i^{(j)}} \pi_{i,k}^{(j)} \mathcal{N}(p|\mu_{i,k}^{(j)}, \sigma_{i,k}^{(j)})
\end{equation}
where $\pi_{i,k}^{(j)}/\mu_{i,k}^{(j)}/\sigma_{i,k}^{(j)}$ 
are the weight/mean/standard deviation of the $k$-th Gaussian component for station $i$ and the $j$-th forecasting interval.
$\mathcal{N}(\cdot)$ denotes the Gaussian distribution.

\subsubsection{Stage 2: Real-time Probabilistic Forecasting}

\textit{\textbf{Step 2.1:}} 
At time slot $t$,
the feature matrix $\bm{X}_t$ is input into the trained deterministic forecasting model,
and the deterministic forecasting result of station $i$ at time slot $t$, i.e., 
$\hat{P}_{i,t}$ is obtained.

\textit{\textbf{Step 2.2:}} 
Check the interval $I^{\mathrm{F}}_j$ that $\hat{P}_{i,t}$ falls into. 
Then, the corresponding forecasting error follows the distribution of $\tilde{\mathcal{P}}_{i}^{(j)}$,
and the probabilistic forecasting of station $i$ at time slot $t$ is given by 
$\hat{\bm{\mathcal{P}}}_{i,t}=\tilde{\mathcal{P}}_{i}^{(j)}+\hat{P}_{i,t}$.

Next, we first analyze the spatio-temporal correlations across multiple EV charging stations, 
and then introduce the detailed structure of the proposed adaptive spatio-temporal graph convolutional network.

\subsection{Spatio-Temporal Correlation Analysis of Multiple EV Charging Stations}
\begin{figure}[htbp]
    \setlength{\abovecaptionskip}{-0.1cm}
    \setlength{\belowcaptionskip}{-0.1cm}
    \centerline{\includegraphics[width=0.95\columnwidth]{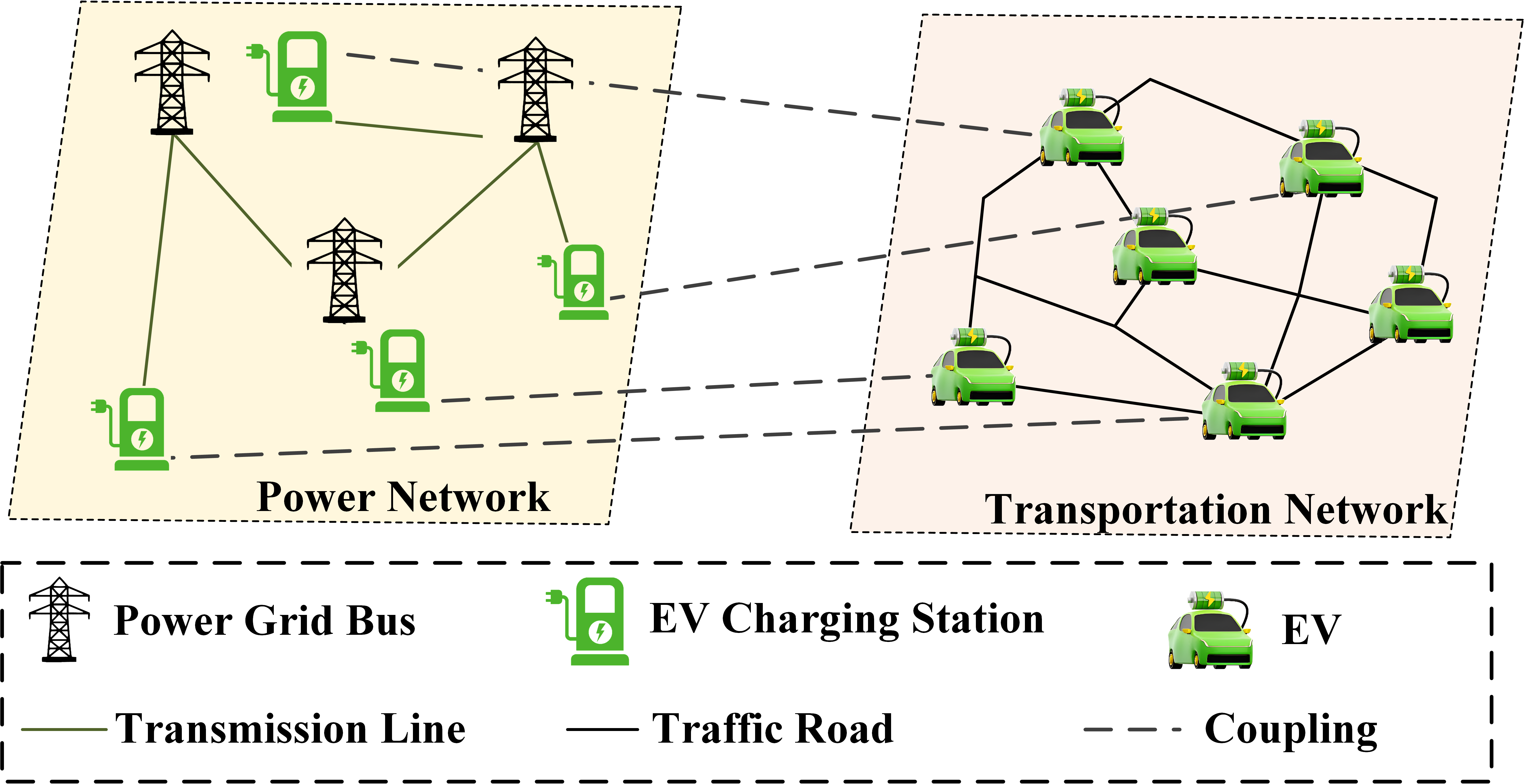}}
    \caption{\rmfamily Coupling between transportation and power networks.}
    \label{DN-EV}
\end{figure}

As transportation electrification advances, 
the integration of an increasing number of EVs into the power system via charging stations intensifies 
the complexity of interactions between transportation and power networks, as illustrated in Fig.~\ref{DN-EV}.
In this context,
the multiple charging stations exhibit complex spatio-temporal correlations.

In the spatial dimension, 
factors within the transportation network, 
such as road topology and infrastructure conditions, 
directly influence EV driving routes and traffic flow, 
which further impact the spatial distribution of potential charging demands in the power network. 
Simultaneously, 
factors in the power network, 
such as the geographical locations of charging stations, 
charging electricity prices, 
and the quality of existing charging facilities, 
also influence the charging choices and driving routes of EVs, 
further affecting traffic flow within the transportation network.

In the temporal dimension, 
real-time road congestion and the availability of charging stations within the transportation network directly impact EV charging schedules, 
influencing traffic flow and subsequently affecting charging demands in the power network. 
In the power network, 
the charging process requires a certain amount of time, 
and factors such as real-time charging prices and the operational status of EVs can influence charging schedules, 
creating cascading effects on the distribution of charging demands.

Furthermore, the influencing factors above can be categorized into time-invariant and time-varying factors.
Time-invariant factors include the geographical locations of the charging stations, the topology of the transportation network, 
and the inherent characteristics of the charging facilities. 
Time-varying factors encompass the real-time traffic flow, weather conditions,
operational status of EVs, and immediate availability of charging facilities. 

In this context, effectively extracting the underlying spatio-temporal features of the multiple EV charging stations 
is essential for precisely forecasting EV charging demands.
 
\subsection{Spatio-Temporal Graph Modeling of Multiple EV Charging Stations}
\begin{figure}[htbp]
    \setlength{\abovecaptionskip}{-0.1cm}
    \setlength{\belowcaptionskip}{-0.1cm}
    \centerline{\includegraphics[width=\columnwidth]{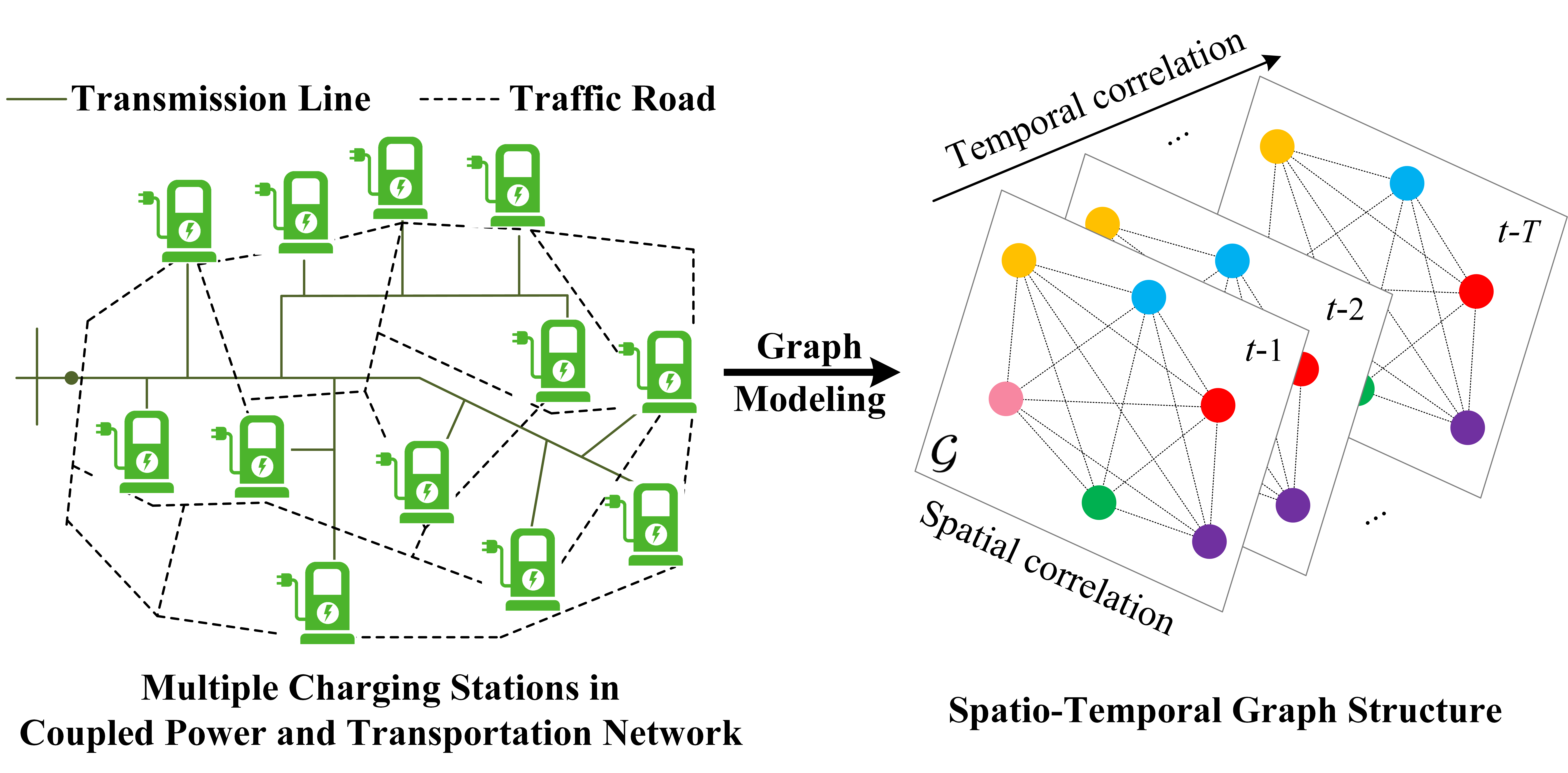}}
    \caption{\rmfamily Spatio-temporal graph modeling of multiple EV charging stations.}
    \label{Graph structure}
\end{figure}

An undirected spatio-temporal graph is adopted 
to model the underlying spatio-temporal coupling features of the multiple charging stations, as illustrated in Fig.~\ref{Graph structure}.
The graph is denoted as $\mathcal{G} = (\mathcal{V}, \mathcal{E})$.
$\mathcal{V}$ is the set of graph nodes corresponding to the charging stations,
and satisfies $|\mathcal{V}| = N^{\mathrm{C}}$.
$\mathcal{E}$ is the set of graph edges representing the underlying spatial-temporal coupling across the charging stations.
Since the charging stations are indirectly connected,
we use the weighted adjacency matrix $\bm{W} \in \mathbb{R}^{N^{\mathrm{C}}\times N^{\mathrm{C}}}$ 
to describe the spatial correlation across the charging stations.

Next, based on the spatio-temporal graph,
we seek to construct an adaptive spatio-temporal graph convolutional network to extract the underlying spatio-temporal features 
of the multiple EV charging stations.
Specifically, 
we propose an adaptive graph convolutional layer to extract spatial features,
which is further combined with an attention-based gated convolutional layer to extract temporal features.
Additionally, a second-order pooling layer is employed to enhance the graph representation.
Finally, the spatio-temporal features are fused to generate the deterministic forecasting results.

\subsection{Adaptive Graph Convolutional Network for Spatial Feature Extraction}
In this subsection, 
we develop an adaptive graph convolutional network to extract spatial features across multiple EV charging stations. 
We begin by introducing an adaptive weighted adjacency matrix that captures both time-invariant and time-varying spatial correlations,
followed by an adaptive graph convolutional layer for efficient feature extraction.

\subsubsection{Adaptive Weighted Adjacency Matrix}
We introduce a learnable adaptive weighted adjacency matrix $\bm{W}$,
which consists of a time-invariant adaptive weighted adjacency matrix $\bm{W}^{\mathrm{TI}}$, 
and a time-varying adaptive weighted adjacency matrix $\bm{W}^{\mathrm{TV}}$, 
modeling the time-invariant and time-varying correlations across multiple EV charging stations, 
respectively. 
\begin{equation}\label{adaptive W}
    \bm{W} = \bm{W}^{\mathrm{TI}} + \bm{W}^{\mathrm{TV}} 
\end{equation}

Compared to the traditional fixed adjacency matrix,
this dual-component structure of $\bm{W}$ provides a more comprehensive and precise representation of the spatial coupling correlations 
by incorporating both time-invariant and time-varying factors. 
Additionally, $\bm{W}$ is optimized end-to-end using stochastic gradient descent, 
which allows the model to adaptively learn the correlation features directly from the data. 
This end-to-end learning enhances the model's ability to identify and capture complex, hidden spatial coupling correlations.

First, leveraging graph embedding theory~\cite{wavenet}, we construct the time-invariant adaptive weighted adjacency matrix $\bm{W}^{\mathrm{TI}}$:
\begin{equation}\label{Ws}
    \bm{W}^{\mathrm{TI}} =  \text{ReLU}( \text{Softmax}( (\bm{E}\bm{E}^T)))
\end{equation}
where $\bm{E} \in \mathbb{R} ^{N^{\mathrm{C}}\times N^{\mathrm{E}}}$ is the embedding matrix of graph $\mathcal{G}$,
and $N^{\mathrm{E}}$ is the embedding dimension.
In terms of activation functions, 
$ \text{ReLU}(\cdot)$ is applied to eliminate weak connections, and $\text{Softmax}(\cdot)$ is applied for weight normalization. 
Formulation~\eqref{Ws} captures the inherent structure and relationships of multiple EV charging stations within the coupling networks
through nonlinear mapping of the multiplied embedding matrix, 
thereby allowing for a nuanced representation that enhances the adaptability.

Second, we construct the time-varying adaptive weighted adjacency matrix $\bm{W}^{\mathrm{TV}}$ 
to model the time-varying coupling across multiple EV charging stations.
Given that the charging demands at the previous time slots provide a realistic reflection of EV travel and charging behavior patterns,
we utilize the generalized similarity of the charging demands across the charging stations to model the time-varying features.
Denoting the charging demands at charging station $i$ and $j$
as $\bm{P}_i=\{P_{i,\tau}|\tau \in [t-T,t-1]\}$ 
and $\bm{P}_j=\{P_{j,\tau}|\tau \in [t-T,t-1]\}$, 
the generalized similarity between the charging demands
is measured using the generalized Mahalanobis distance:
\begin{equation}
    d^\mathrm{M} (\bm{P}_i, \bm{P}_j) 
    = \sqrt{(\bm{P}_i - \bm{P}_j) \bm{M}\bm{M}^\text{T} (\bm{P}_i - \bm{P}_j)^\text{T}}
\end{equation}
where $\bm{M} \in \mathbb{R}^{T \times T} $ is a learnable matrix that captures the temporal correlations across the multiple EV charging stations.
Then, the time-varying adaptive weighted adjacency matrix $\bm{W}^{\mathrm{TV}}$ is densified
by applying the normalized Gaussian kernel function:
\begin{equation}
    \bm{W}^{\mathrm{TV}}\! =\! \left\{ W^{\mathrm{TV}}_{ij} \middle| W^{\mathrm{TV}}_{ij}\! = \!
    \frac{\exp\left(-d^\mathrm{M} (\bm{P}_i, \bm{P}_j)/(2\sigma^2)\right)}
    {\sum\limits_{k \in \mathcal{V}} \exp\left(-d^\mathrm{M} (\bm{P}_i, \bm{P}_k)/(2\sigma^2)\right)},\! \ i,j\! \in\! \mathcal{V} \right\}
\end{equation}
where $\sigma$ is the bandwidth of the Gaussian kernel function,
and $\exp(\cdot)$ denotes the exponential function.

\subsubsection{Spectral Graph Convolutional Network}
Graph convolutional networks (GCNs) are widely used for extracting spatial features from graph-structured data. 
GCNs are typically categorized into two main types: 
spatial domain and spectral domain. 
Spatial-based methods define graph convolutions by aggregating features from neighboring nodes, 
whereas spectral-based methods use filters based on graph signal processing principles. 
By performing convolutions in the Fourier domain, 
spectral GCNs have the advantage of capturing global information from the graph and offering relatively easy computation. 
In this paper, 
we employ spectral GCN to leverage these benefits.
The core formula is as follows:
\begin{equation}\label{GCN0}
    \bm{X}*_{\mathcal{G}}\bm{g}_{\theta} = \sum_{{k}=1}^{K^\mathrm{s}} \theta_{k} T_{k}(\tilde{\bm{L}}) \bm{X}
\end{equation}
where $\tilde{\bm{L}} = 2\bm{L} / \lambda^{\mathrm{max}} - \bm{I}$ is the normalized adaptive graph Laplacian matrix.
$\bm{L} = \bm{I} - \bm{(D)^{-\frac{1}{2}}} \bm{W} \bm{(D)^{-\frac{1}{2}}}$.
$\bm{D}= \{D_{ij} = \sum_j W_{ij} |i,j\in\mathcal{V}\}$ is the degree matrix of the graph,
and $\bm{I}\in \mathbb{R}^{N^{\mathrm{C}}\times N^{\mathrm{C}}}$ is the identity matrix.
$\lambda^{\mathrm{max}}$ is the largest eigenvalue of $\bm{L}$.
$\bm{X}$ is the input graph feature matrix.
$g_{\theta}$ is the filter parameterized by $\theta$.
$\theta_k$ is the Chebyshev coefficients.
$T_{k}(\cdot)$ is the Chebyshev polynomial of order $k$.
$K^\mathrm{s}$ is the number of Chebyshev polynomials.
The details of the spectral GCN can be referred in~\cite{spectralGCN,wu2020comprehensive}.

The spatial feature extraction can be conducted by~\eqref{GCN0}.
Next, we introduce the attention-based gated convolutional network to extract temporal features.

\subsection{Attention-based Gated Convolutional Neural Networks for Temporal Feature Extraction}
\subsubsection{Gated Convolutional Neural Networks}
Recurrent neural networks have become a cornerstone in time series analysis due to their ability to capture sequential dependencies. 
However, traditional recurrent neural networks models, including long short-term memory networks and bidirectional LSTM networks, 
face significant challenges from the inherently sequential nature of their iteration process, 
which are complex and time-consuming.
To address these limitations, we draw inspiration from the advancements in~\cite{Gated_CNN} 
and propose the usage of gated convolutional neural networks and gated linear units (GLUs) along the temporal axis to extract temporal features.
The formulation is as follows:
\begin{equation}\label{GCNN}
    \bm{X}^{\mathrm{out}} = (\bm{X}^{\mathrm{in}} \ast \bm{B} + \bm{b}) \odot \psi (\bm{X}^{\mathrm{in}} \ast \bm{C} + \bm{c})
\end{equation}
where $\bm{X}^{\mathrm{in}} \in \mathbb{R}^{C^{\mathrm{in}} \times N^{\mathrm{C}} \times T}$ is the input feature matrix,
$\bm{X}^{\mathrm{out}}$ is the output feature matrix,
and $C^{\mathrm{in}}/C^{\mathrm{out}}$ are the number of input/output channels.
$\bm{B}/\bm{C}$ and $\bm{b}/\bm{c}$ are the learnable convolutional kernels and biases.
$\psi$ is the activation function.

Eq.\eqref{GCNN} leverages the strengths of convolutional operations to capture temporal features more efficiently than traditional RNNs. 
GLUs can effectively address the vanishing gradient problem that often plagues deep neural networks by providing a linear path for gradient flow, 
while still maintaining the non-linear capabilities necessary for complex feature representation. 
This combination of techniques allows our model to efficiently learn temporal dependencies in time series. 
Furthermore, 
by incorporating residual connections into the network, 
we aim to mitigate the risk of overfitting and enhance the model's ability to generalize to unseen data samples.

\subsubsection{Temporal Attention Mechanism}

In EV charging demands time series, 
each time slot reflects unique patterns and trends of EV usage. 
Attention mechanisms are used to capture the varying importance of time slots by assigning different weights. 
This ensures that the model focuses on periods most indicative of future trends.
Traditional attention mechanisms face a significant limitation of 
quadratically increasing computational complexity with the length of the sequence,
which can be a bottleneck for long time series.
To address this challenge, 
we employ the scaled dot-product attention mechanism:
\begin{equation}\label{QKVattention}
    \bm{A^{\mathrm{tt}}} = \text{Softmax}\left(\frac{\bm{Q}\bm{K}^T}{\sqrt{d_k}}\right)\bm{V}
\end{equation}

Scaled dot product attention mechanism takes $\bm{Q}$, $\bm{K}$ and $\bm{V}$ as inputs, which represent the queries, keys, and values, respectively. 
$d_k$ is the dimension of queries and keys.
Dot products between keys and queries are computed before being passed to a Softmax function to obtain the attention weights, 
which are then multiplied by the values to produce the final outputs.
This mechanism can effectively capture the differential importance of various time slots in the EV charging demands time series,
enhancing the model's ability to accurately capture temporal dependencies while also maintains computational efficiency, 
making it feasible to handle extensive sequences of time series.

\subsection{Spatio-Temporal Graph Convolutional Block}
\begin{figure}[htbp]
    \setlength{\abovecaptionskip}{-0.1cm}
    \setlength{\belowcaptionskip}{-0.1cm}
    \centerline{\includegraphics[width=0.95\columnwidth]{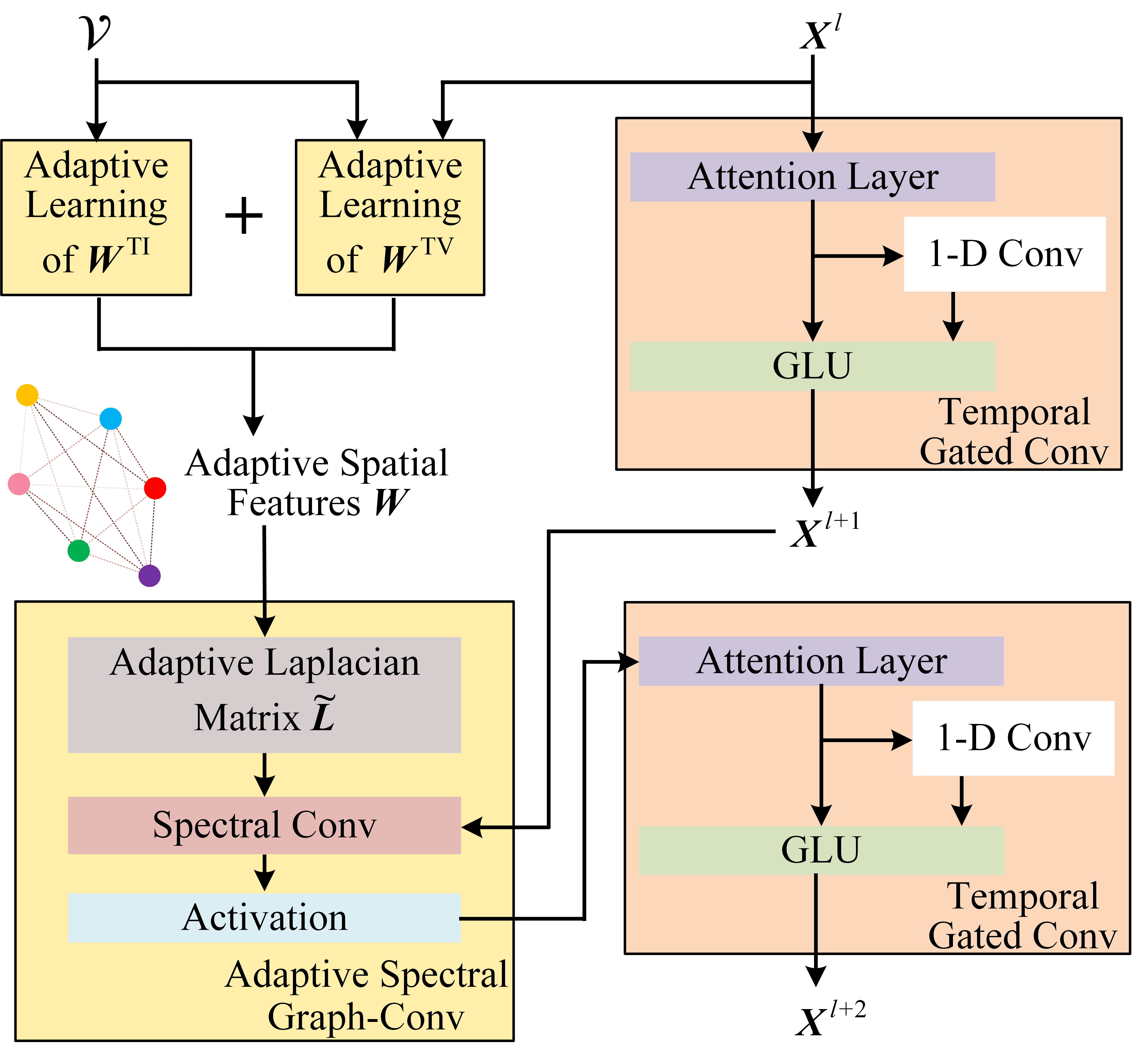}}
    \caption{\rmfamily Structure of the spatio-temporal graph convolutional block.}
    \label{ST-Conv}
\end{figure}

To simultaneously extract spatial and temporal features, 
a spatio-temporal graph convolutional (ST-Conv) block is constructed by 
combining adaptive spectral graph convolutional layers with attention-based temporal gated convolutional layers. 
The ST-Conv is structured by a graph layer positioned between two temporal layers, 
facilitating efficient spatial information dissemination across temporal convolutions. 
This design enables jointly extracting spatial and temporal dependencies. 
The schematic representation of the ST-Conv block is illustrated in Fig.~\ref{ST-Conv}.

Within the spatio-temporal graph convolutional block, 
the input feature matrix $\bm{X}^l$ undergoes temporal feature extraction through the first temporal layer,
and generates the intermediate feature matrix $\bm{X}^{l+1}$.
Subsequently, 
$\mathcal{V}$ is combined with $\bm{X}^l$ for adaptive learning of $\bm{W}^{\mathrm{TI}}$ and $\bm{W}^{\mathrm{TV}}$. 
In the adaptive spectral graph convolutional layer, 
the input feature matrix $\bm{X}^{l+1}$ is fused with the learnable weighted adjacency matrix $\bm{W}$ 
to extract spatial features. 
Following activation, the processed data traverses another temporal convolutional layer, 
ultimately yielding the extraction of integrated spatio-temporal features.

\subsection{Second-Order Pooling for Graph Presentation}
Graph pooling aggregates node features within a graph to generate a unified representation of the entire graph structure.
Traditional graph pooling methods, such as max pooling, average pooling, and sum pooling, 
only collect first-order statistics and disregard feature correlation information, 
which may result in inadequate handling of variable node representations and consequently hinder overall model performance. 
In this paper, 
we introduce a bilinear mapping second-order pooling method~\cite{SOP} to enhance feature extraction 
by capturing second-order feature correlations and topological information from all nodes while reducing the output dimensionality. 
In the bilinear mapping second-order pooling layer,
the graph representation vector $h_{\mathcal{G}}$ is computed through:
\begin{equation}
    h_{\mathcal{G}} = \text{flatten}(\bm{Z}^T\bm{X}^T\bm{X}\bm{Z})
\end{equation}

\begin{figure}[htbp]
    \setlength{\abovecaptionskip}{-0.1cm}
    \setlength{\belowcaptionskip}{-0.1cm}
    \centerline{\includegraphics[width=0.95\columnwidth]{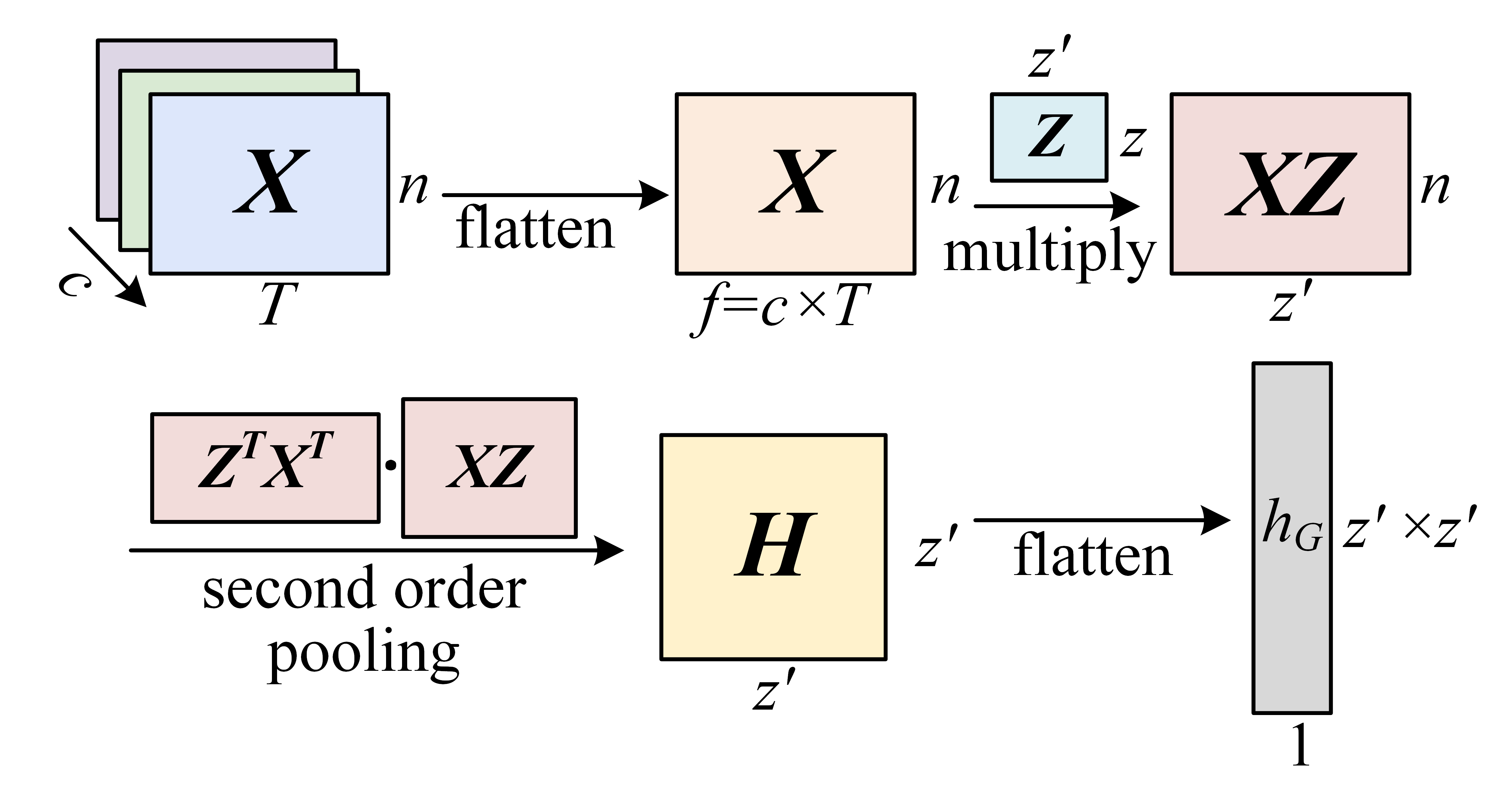}}
    \caption{\rmfamily Second-order pooling for graph representation.}
    \label{SOP}
\end{figure}

Here, $\text{flatten}(\cdot)$ is a function that reshapes the matrix into a vector.
$\bm{X}$ represents the graph feature matrix after applying ST-Conv blocks for feature extraction, 
and $\bm{Z}\in \mathbb{R}^{z'\times z}$ is the linear mapping matrix. 
The procedure for second-order pooling is illustrated in Fig.~\ref{SOP}.

\subsection{Output Layer}
Based on $h_{\mathcal{G}}$, 
the output layer employs a multi-layer perceptron (MLP) block composed of linear mappings, batch normalization, 
and activation functions. 
$\mathrm{MLP}(\cdot)$ extracts relevant features to generate deterministic forecasting of multiple EV charging demands:
\begin{equation}
    [\hat{P}_{i,t}]_{i\in \varOmega^{\mathrm{C}}} = \mathrm{MLP}(h_{\mathcal{G}})
\end{equation}

\subsection{Overall Model Structure of ASTGCN}
\begin{figure}[htbp]
    \setlength{\abovecaptionskip}{-0.1cm}
    \setlength{\belowcaptionskip}{-0.1cm}
    \centerline{\includegraphics[width=0.95\columnwidth]{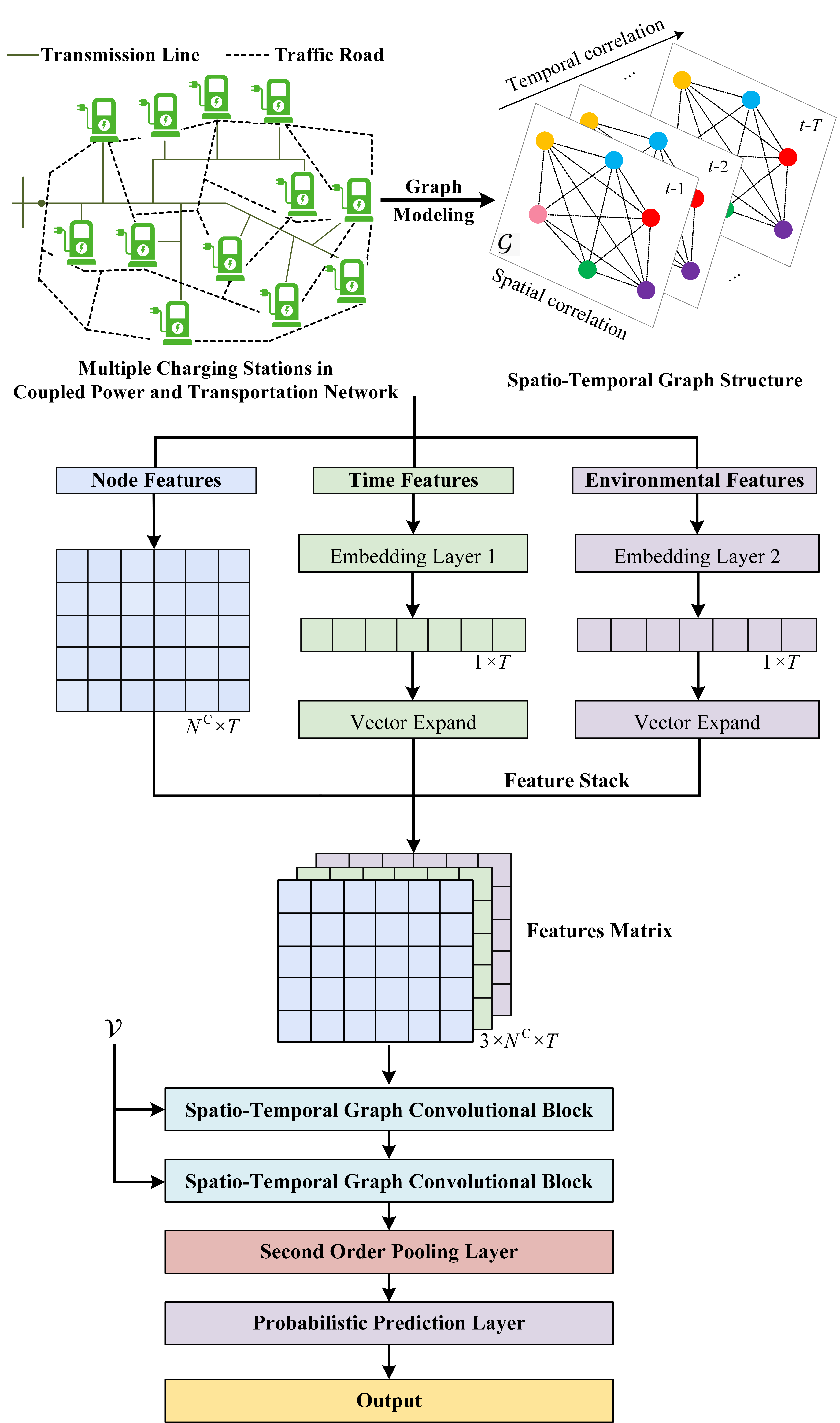}}
    \caption{\rmfamily Diagram of the ASTGCN model structure.}
    \label{ASTGCN} 
\end{figure}

By integrating the aforementioned feature extraction layers, 
we develop an adaptive spatio-temporal graph convolutional neural network with attention mechanism
for deterministic forecasting of multiple EV charging demands.
The overall structure of ASTGCN is shown in Fig.~\ref{ASTGCN}.

The input temporal features primarily consist of node features (e.g., historical observations of charging demands), 
time features (e.g., time slots and days of the week), and environmental features (e.g., weather conditions). 
The covariates can contribute to improving the accuracy of probabilistic forecasting 
by capturing the underlying patterns and relationships in the data.
The covariates can be selected based on the availability of data and the characteristics of the forecasting task.

Initially, the time and environmental features are processed through an embedding layer to generate feature embedding vectors. 
These vectors are then reshaped and expanded to match the shape of the node features. 
Subsequently, the node, time, 
and environmental features are integrated to construct the final feature matrix.
The integrated features are first fed into two spatio-temporal graph convolutional blocks to extract both spatial and temporal features. 
In the ST-Conv block, spatial features are captured by the adaptive graph convolutional layer, 
while temporal features are extracted using the attention-based gated convolutional layer. 
The extracted spatio-temporal features are then processed by a second-order pooling layer, 
generating graph representation by capturing second-order information.
Finally, this graph representation is fed into the output layer to produce probabilistic forecasts of the EV charging demands.
 
\subsection{Loss Function}

The loss function of ASTGCN
is defined as the mean squared error (MSE) between the forecasting results and the ground truth:
\begin{equation}
    \mathrm{Loss} = \frac{1}{|\mathcal{D}|}\sum_{P_{i,t} \in \mathcal{D}}\sum_{i\in \varOmega^{\mathrm{C}}} 
     (P_{i,t} - \hat{P}_{i,t})^2
\end{equation}
where $\mathcal{D}$ is the training dataset with size $|\mathcal{D}|$.
The summation over $P_{i,t} \in \mathcal{D}$ spans all data samples in the dataset, 
over $i$ spans all charging stations. 

The training process of ASTGCN is conducted by minimizing the loss function using stochastic gradient descent,
in which process the parameters of the model are updated and optimized.

\section{Real-Time Risk Analysis Based on Probabilistic Power Flow}\label{Probabilistic Risk Analysis}
The probabilistic characteristics of EV charging demands lead to probabilistic power flow within the distribution network. 
Based on the probabilistic forecasting of EV charging demands, 
this section presents an improved PPF calculation method for real-time risk analysis in distribution networks 
with large-scale EV integration.
Note that we focus on investigating the initial behavior and corresponding inherent risks of EVs
without considering flexible charging or discharging of EVs~\cite{EVoperation1,EVoperation2},
or the coordination with other flexible resources~\cite{EVAggregate,EVoperation3,EVoperation4}.
In this context, we assume that the regular load demands at other buses are constant and predetermined.
Among the risk criterions,
the low voltage criterion is a significant concern due to the heavy load imposed on the power system by massive EV integration 
and the strict regulations enforced for low voltage conditions~\cite{HCTL}.

\subsection{GMM-Based Probabilistic Power Flow}
The power flow equations can be described as:
\begin{equation}\label{pf}
    \bm{w} = f(\bm{x})
\end{equation}
where $\bm{w}$ represents the active and reactive power injections of buses, 
$\bm{x}$ represents the bus voltage angles and magnitudes, 
and $f(\cdot)$ represents the power flow equations.
Taylor series expansion is carried out at the base operating point $\bm{x}^{\mathrm{0}}$, 
and higher-order terms of 2 or more are ignored.
In this case,~\eqref{pf} can be simplified using the Jacobian matrix~\cite{linear_pf1,linear_pf2}:
\begin{subequations}
    \begin{align}
        \label{pf1}
        &\Delta \bm{w} = -\bm{J}\Delta \bm{x}\\
        \label{pf2}
        &\Delta \bm{x} = (-\bm{J})^{-1}\Delta \bm{w} = \bm{S}\Delta \bm{w}
    \end{align}
\end{subequations}
where $\bm{J}/\bm{S}$ are the Jacobian matrix/sensitivity matrix of the power flow equations
at the last iteration.
$\Delta \bm{w}$ represents the random power injection disturbances of buses, 
and $\Delta \bm{x}$ represents the resulting variations in bus voltages.

Recall that the EV charging demands are forecasted as GMMs in~\eqref{GMMforecast},
combining~\eqref{pf2} and~\eqref{GMMforecast}, we have:
\begin{equation}\label{BGMMpf}
    \Delta x_i \!\sim \!\mathcal{X}_i\!=\!\sum_{\mathbf{k}}\! \left( \left( \prod_{s=1}^{|\bm{S}|} \pi_{s k_s}^{(j_s)} \right) 
    \mathcal{N}\! \left( \sum_{s=1}^{|\bm{S}|} \bm{S}_{is} \mu_{s k_s}^{(j_s)}, 
    \sqrt{\sum_{s=1}^{|\bm{S}|} \bm{S}_{is}^2 (\sigma_{s k_s}^{(j_s)})^2}\! \right) \!\right)
\end{equation}
where $\mathbf{k} = (k_1, k_2, \dots, k_{|\bm{S}|})$ 
is a combination of realizations from $[K_1] \times \dots \times [K_{|\bm{S}|}]$.
$|\bm{S}|$ denotes the dimension of $\bm{S}$.
Notably, $\pi_{s k_s}^{(j_s)},\mu_{s k_s}^{(j_s)},\sigma_{s k_s}^{(j_s)}=0$ if $s \notin \varOmega^{\mathrm{C}}$.
\eqref{BGMMpf} provides the analytical expression of the probabilistic power flow.
However, the number of components in~\eqref{BGMMpf} is $\prod_{s=1}^{|\bm{S}|} K_s$, 
which makes~\eqref{BGMMpf} computationally intensive and potentially introduce scalability challenges. 
Some studies have proposed to reduce the number of components in the GMMs to enhance computational efficiency while preserving accuracy.
Technic details can be referred in~\cite{GMMreduction1,GMMreduction2}.

\subsection{Operational Boundaries Identification}
Based on~\eqref{BGMMpf}, 
the voltage $\bm{x} $ follows $\bm{x}=\bm{x}^{\mathrm{0}} + \Delta \bm{x}$.
Next,
we seek to identify the operation voltage lower boundary $\underline{V}$.
First, we define a chance constraint voltage constraint violation limit as $\varsigma$,
which is defined as:
\begin{equation}\label{operational boundaries}
    \mathcal{X}_i(x_i \leq \underline{V}) \leq \varsigma, \forall i \in \varOmega^{\mathrm{L}}
\end{equation}
where $\varOmega^{\mathrm{L}}$ is the set of load buses.
$\underline{V}$ can be identified by solving~\eqref{operational boundaries}.

\section{Real-time HC Assessment for EVs}\label{Real-time Hosting Capacity Assessment for EVs}
In this section, we propose an optimization model to assess the real-time HC of EVs within the distribution network
considering the real-time stochasticity of EV charging demands.
\subsection{Model Formulation}

\subsubsection{Objective Function}
The traditional HC of EVs assessed from a long-term perspective is defined as 
the maximum EV charging demands that the distribution network can accommodate without violating network operation constraints,
as in~\eqref{HC1}.
\begin{equation}\label{HC1}
\begin{split}
    &\max \sum_{i\in \varOmega^{\mathrm{C}}} \overline{P}_{i}\\
    &s.t., \text{network~operation~constraints}
\end{split}
\end{equation}
where $\overline{P}_{i}^{\mathrm{EV}}$ is the upper limit of accommodated EV charging demands at charging station $i$.

In this paper, 
we focus on the real-time stochasticity of EV charging demands.
By modeling them as stochastic variables through probabilistic forecasting,
we define the real-time HC of EVs at time slot $t$
as the maximum expected EV charging demands satisfaction that the distribution network can accommodate 
without violating network secure operation constraints.
The objective of real-time HC assessment is:
\begin{equation}\label{HC2}
    \max \sum_{i\in \varOmega^{\mathrm{C}}} \int_{0}^{\overline{P}_{i,t}}p_{i,t}\hat{\mathcal{P}}_{i,t}(p_{i,t})dp_{i,t}
\end{equation}
where $\hat{\mathcal{P}}_{i,t}$ is the forecasted probabilistic distribution of EV charging demands at charging station $i$ 
and time slot $t$.

Recall that we focus on investigating the initial inherent risks of EVs
without considering flexible operations of EVs or any other flexible resources.
The network operation constraints mainly include the power flow constraints 
and the minimum EV charging demands satisfaction constraints.
\subsubsection{Power Flow Constraints}
DistFlow branch model is employed to model the power flow of distribution network:
\begin{subequations}\label{opf}
    \begin{align}
        \label{pf U}
        & V_{j,t}^2 \!=\! V_{i,t}^2 \!-\! 2(r_{ij}P_{ij,t}^{\mathrm{Ln}} \!+\! x_{ij}Q_{ij,t}^{\mathrm{Ln}})\!+\!(r_{ij}^2 + x_{ij}^2)I_{ij,t}^2 \\
        \label{pf P}
        & p_{j,t} = P_{ij,t}^{\mathrm{Ln}}-r_{ij}I_{ij,t}^2-\sum_{l:j\rightarrow l}P_{jl,t}^{\mathrm{Ln}} \\
        \label{pf Q}
        & q_{j,t} = Q_{ij,t}^{\mathrm{Ln}}-x_{ij}I_{ij,t}^2-\sum_{l:j\rightarrow l}Q_{jl,t}^{\mathrm{Ln}} \\
        \label{pf nonconvex}
        & V_{i,t}^2I_{ij,t}^2 = (P_{ij,t}^{\mathrm{Ln}})^2+(Q_{ij,t}^{\mathrm{Ln}})^2 \\
        \label{pf U range}
        & \underline{V}\leq V_{i,t}\leq \overline{V} \\
        \label{pf I range}
        & |I_{ij,t}| \leq\overline{I}_{ij} \\
        \label{pin EV}
        & p_{i,t} = P_{i,t}^{\mathrm{L}}+\overline{P}_{i,t}, \forall i \in \varOmega^{\mathrm{C}}\\
        \label{pin}
        & p_{i,t} = P_{i,t}^{\mathrm{L}}, \forall i \in \varOmega^{\mathrm{L}}/\varOmega^{\mathrm{C}}\\
        \label{qin}
        & q_{i,t}= Q_{i,t}^{\mathrm{L}} 
    \end{align}
\end{subequations}

$r_{ij}/x_{ij}$ are the line resistance/reactance of line $ij$, respectively.
$I_{ij}$ is the electric current of line $ij$, 
with $\overline{I}_{ij}$ as the upper line current.
$V_{i}$ is the voltage of bus $i$, with $\overline{V}/\underline{V}$ as the upper/lower bus voltage.
$P_{ij}^{\mathrm{Ln}}/Q_{ij}^{\mathrm{Ln}}$ are the line active/reactive power of line $ij$, respectively.
$p_{i}/q_{i}$ are the active/reactive injection power of bus $i$.
$P_{i}^{\mathrm{L}}/Q_{i}^{\mathrm{L}}$ are the active/reactive load power of bus $i$.
$\varOmega^{\mathrm{L}}/\varOmega^{\mathrm{C}}$ is the set of load buses excluding the buses with charging stations.
\eqref{pf U} describes the voltage drop over line $ij$. 
\eqref{pf P} and~\eqref{pf Q} represent the active and reactive power balance of bus $j$.
\eqref{pf nonconvex} is the power flow equation of line $ij$.
\eqref{pf U range},~\eqref{pf I range} are the security constraints 
and we mainly concentrate on the voltage security constraint in this paper.
\eqref{pin},~\eqref{pin EV} and~\eqref{qin} are the power balance constraints.

\subsubsection{EV Charging Demands Satisfaction Constraints}
To ensure a certain satisfaction level of EV charging demands, 
the maximum accommodated EV charging demands must satisfy specific probabilistic constraints:
\begin{align}\label{satisfy1}
    \hat{\mathcal{P}}_{i,t}(p_{i,t}\leq \overline{P}_{i,t})\geq 1-\epsilon_i,\forall i \in \varOmega^{\mathrm{C}}
\end{align}
where $\epsilon_i$ is the permissible probability of charging demands not being satisfied at bus $i$.
\eqref{satisfy1} can be transformed into a more tractable formulation as:
\begin{align}\label{satisfy2}
    \overline{P}_{i,t}\geq \hat{\mathcal{F}}_{i,t}^{-1}(1-\epsilon_i),\ \forall i \in \varOmega^{\mathrm{C}}
\end{align}
where $\hat{\mathcal{F}}_{i,t}$ is the cumulative distribution function of $\hat{\mathcal{P}}_{i,t}$,
and $\hat{\mathcal{F}}_{i,t}^{-1}$ is the inverse function of $\hat{\mathcal{F}}_{i,t}$.

\subsection{Reformulation of Non-convex Constraints }
The expressions in \eqref{HC2} and \eqref{pf nonconvex} are non-convex and are challenging to be solved directly. 
To tackle this, we reformulate them into more tractable forms.
\subsubsection{Integration in the Objective Function}
Given that each integration term in objective function~\eqref{HC2} increases with the rising value of 
$\overline{P}_{i,t}$,
we can approximate each integration term with a piecewise linear function.

Generally, a integration term $\int_{\underline{p}}^{\overline{p}}p{\mathcal{P}}(p)dp$
can be approximated by a $n$ segment piecewise linear function as:
\begin{equation}\label{piecewise1}
    \int_{\underline{p}}^{\overline{p}}p{\mathcal{P}}(p)dp\approx  f(p)=
    \left\{\begin{matrix}
        k_{1}p+b_{1},& p_{1}\leq p\leq p_{2}  \\
        k_{2}p+b_{2},& p_{2}\leq p\leq p_{3}\\
        \dots ,& \dots\\
        k_{n}p+b_{n}, & p_{n}\leq p\leq p_{n+1}
    \end{matrix}\right.
\end{equation}
where $k_{i}$ and $b_{i}$ are the slope and intercept of the $i$-th segment of the piecewise linear function.
$\underline{p}=p_{1}$ and $\overline{p}=p_{n+1}$ are the lower and upper limits of integration variable $p$, respectively.
Then, by introducing nonnegative continuous variable $\{w_i\}_{i=1}^{n+1} $ and the binary variable $\{z_i\}_{i=1}^{n}$,
$ f(p) $ is represented as:
\begin{subequations}\label{piecewise2}
    \begin{align}
    &f(p) = \sum_{k=1}^{n+1} w_k f(p_k),p = \sum_{k=1}^{n+1} w_k p_k \\
    &w_1 \leq z_1, w_{n+1} \leq z_{n} \\
    &w_k \leq z_{k-1} + z_k,\  \forall k\in [2,n] \\
    &\sum_{k=1}^{n+1} w_k = 1 , \sum_{k=1}^{n} z_k = 1
\end{align}
\end{subequations}

In this way, the non-convex integration term can be convexified to a mixed-integer linear formulation.

\subsubsection{Nonconvex Power Flow Constraint}

The nonconvex constraint~\eqref{pf nonconvex} can be relaxed to the following second-order cone formulation 
by introducing two slack variables $\mathsf{V}_{i}$ and $\mathsf{I}_{ij}$ to replace the quadratic term~\cite{SOCP1,SOCP2,SOCP3}, 
as in~\eqref{SOCP}. 
\begin{subequations}\label{SOCP}
    \begin{align}
        &\left \| \begin{matrix}2P_{ij}^{\mathrm{Ln}}
            \\2Q_{ij}^{\mathrm{Ln}}
            \\\mathsf{V}_{i} - \mathsf{I}_{ij}
        \end{matrix} \right \| _2 \leq \mathsf{V}_{i} + \mathsf{I}_{ij}\\
        &\mathsf{V}_{i} = V_{i}^2,\ \mathsf{I}_{ij} = I_{ij}^2
    \end{align}
\end{subequations}

\subsection{Overall Model}

Finally, the problem of real-time HC assessment for EVs can be formulated as a mixed-integer second-order cone programming problem in~\eqref{opt}, 
which can be efficiently solved by commercial solvers.
\begin{equation}\label{opt}
    \begin{aligned}
        &\max \sum_{i \in \varOmega^{\mathrm{C}}} f_i(p_{i,t})\\
        &\text{s.t.} \ \eqref{pf U}-\eqref{pf Q},\eqref{pf U range}-\eqref{qin}, \eqref{satisfy2}-\eqref{SOCP}
    \end{aligned}
\end{equation}
where $f_i(\cdot)$ is the approximated piecewise linear function of 
$\int_{0}^{\overline{P}_{i,t}}p_{i,t}\hat{\mathcal{P}}_{i,t}(p_{i,t})dp_{i,t}$,
and can be obtained by~\eqref{piecewise1} and~\eqref{piecewise2}.
In~\eqref{opt}, variables include
$\overline{P}_{i,t}$ (i.e., the maximum accommodated EV charging demands at bus $i$),
variables related to power flow equations in~\eqref{opf},
and variables related to piecewise linear function in~\eqref{piecewise1} and~\eqref{piecewise2}.
Parameters include $\{\hat{\mathcal{P}}_{i,t}\}_{i \in \varOmega^{\mathrm{C}}}$ obtained from probabilistic forecasting,
$\{\epsilon_i\}_{i \in \varOmega^{\mathrm{C}}}$, 
parameters related to power flow equations in~\eqref{opf},
and parameters related to piecewise linear function in~\eqref{piecewise1} and~\eqref{piecewise2}.

\section{Numerical Experiments}\label{case study}
In this section, we conduct numerical experiments to evaluate the performance of the proposed 
real-time HC assessment method for EVs.
\subsection{Set up}
\begin{figure}[htbp]
    \setlength{\abovecaptionskip}{-0.1cm}
    \setlength{\belowcaptionskip}{-0.1cm}
    \centerline{\includegraphics[width=0.95\columnwidth]{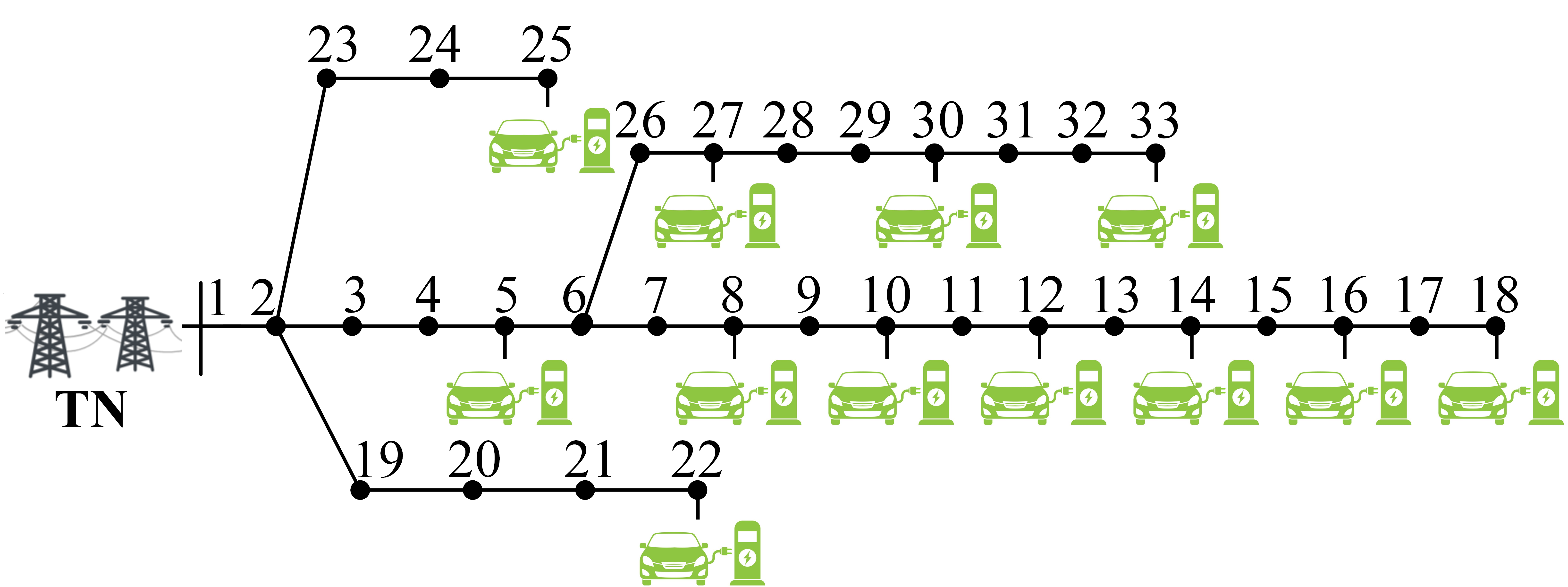}}
    \caption{\rmfamily Structure of the distribution network with multiple charging stations.}
    \label{DNEV}
\end{figure}

We utilize a real-world EV charging dataset of 12 charging stations from Tianjin, China, 
encompassing historical charging transactions throughout the entire year of 2023.
The 12 charging stations are set to be connected to the IEEE 33-bus distribution network
in buses $\varOmega^{\mathrm{C}}=$\{5, 8, 10, 12, 14, 16, 18, 22, 25, 27, 30, 33\},
as shown in Fig.~\ref{DNEV}.

We first process the raw transaction data to generate the charging demands time series for each station
with a 15-minute time interval.
In this paper, we forecast for the next time slot using the observation of past 2 hours,
i.e., $T=8$.
Data cleaning procedures are implemented to eliminate invalid transactions, 
specifically those with energy consumption less than 1 kWh, 
charging durations shorter than 1 minute or longer than 24 hours, 
and transactions exhibiting abnormal charging power. 
Additionally, the power series in the dataset are normalized within $\left[0,1\right]$.
We set $N^{\mathrm{F}}=100$ with interval length of $0.01$.
Then, we randomly partition the dataset into three subsets: 60\% for training, 20\% for validation, and 20\% for testing. 
The models are implemented in Python 3.10.11 and PyTorch 2.3.0 using the RMSprop optimizer and CUDA 12.1 libraries. 
The experiments are run on a server with an NVIDIA GeForce RTX 3070 GPU, 
an Intel Core i9-13900K CPU (5.4GHz), 64 GB of RAM, 
and Windows 11 64-bit operating system.
The training hyperparameters are set as: learning rate $0.001$, batch size $64$, epochs $60$, kernel size $K^\mathrm{s}=3$.
The model achieving the lowest loss on the validation set is saved as the optimal model.

\subsection{Validation Metrics of Deterministic Forecasting}
For deterministic forecasting, we adopt mean absolute error (MAE), root mean squared error (RMSE),
and weighted absolute percentage error (WAPE)~\cite{prob_pred_PV} as the validation metrics. 
Unlike the mean absolute percentage error (MAPE), 
WAPE offers greater robustness when true values approach zero, 
as it uses the sum of true values in the denominator 
thereby mitigating the influence of extreme values and providing a more stable measure of error.
\begin{subequations}
    \begin{align}
        &\mathrm{MAE} = \frac{1}{|\mathcal{D}|}\sum_{P_{i,t} \in \mathcal{D}}\sum_{i\in \varOmega^{\mathrm{C}}} 
        |P_{i,t} - \hat{P}_{i,t}|\\
        &\mathrm{RMSE} = \sqrt{\frac{1}{|\mathcal{D}|}\sum_{P_{i,t} \in \mathcal{D}}\sum_{i\in \varOmega^{\mathrm{C}}} 
        (P_{i,t} - \hat{P}_{i,t})^2}\\
        &\mathrm{WAPE} = \frac{\sum_{P_{i,t} \in \mathcal{D}}\sum_{i\in \varOmega^{\mathrm{C}}}|P_{i,t} - \hat{P}_{i,t}|}
                     {\sum_{P_{i,t} \in \mathcal{D}}\sum_{i\in \varOmega^{\mathrm{C}}} 
                     P_{i,t} }
    \end{align}
\end{subequations}

\subsection{Deterministic Forecasting Performance of ASTGCN}
\begin{figure}[htbp]
    \setlength{\abovecaptionskip}{-0.1cm}
    \setlength{\belowcaptionskip}{-0.1cm}
    \centerline{\includegraphics[width=0.95\columnwidth]{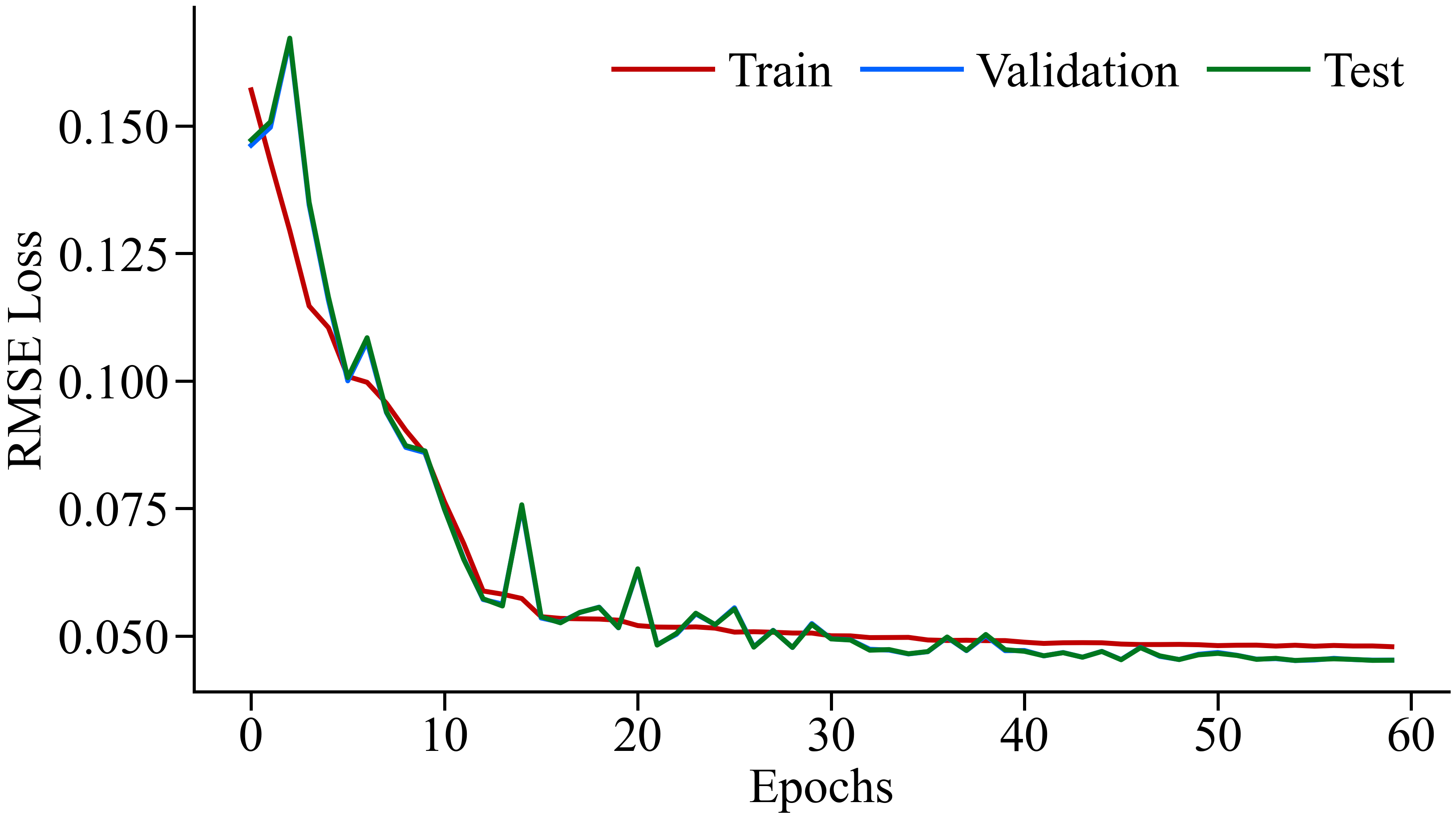}}
    \caption{\rmfamily Loss curves on training, validation, and test sets.}
    \label{predloss}
\end{figure}

The loss curves on the training, validation, and test sets during model training are shown in Fig.~\ref{predloss}.
The curves show a rapid loss decline in the initial epochs, 
followed by a gradual decrease to a stable convergence as training progresses. 
This suggests that the model quickly learns the data patterns and effectively reduces forecasting error in the early training stages. 
In the later stages, 
the model maintains stable performance without significant overfitting or underfitting. 
Furthermore, 
the close alignment of validation and test losses with the training loss demonstrates the model's strong generalization ability.
Finally, the proposed ASTGCN model achieves a low RMSE loss of 0.0442 on the test set.

\begin{figure}[htbp]
    \centering
    \begin{subfigure}{\columnwidth}
        \setlength{\abovecaptionskip}{0.1cm}
        \setlength{\belowcaptionskip}{0.1cm}
        \centerline{\includegraphics[width=0.95\columnwidth]{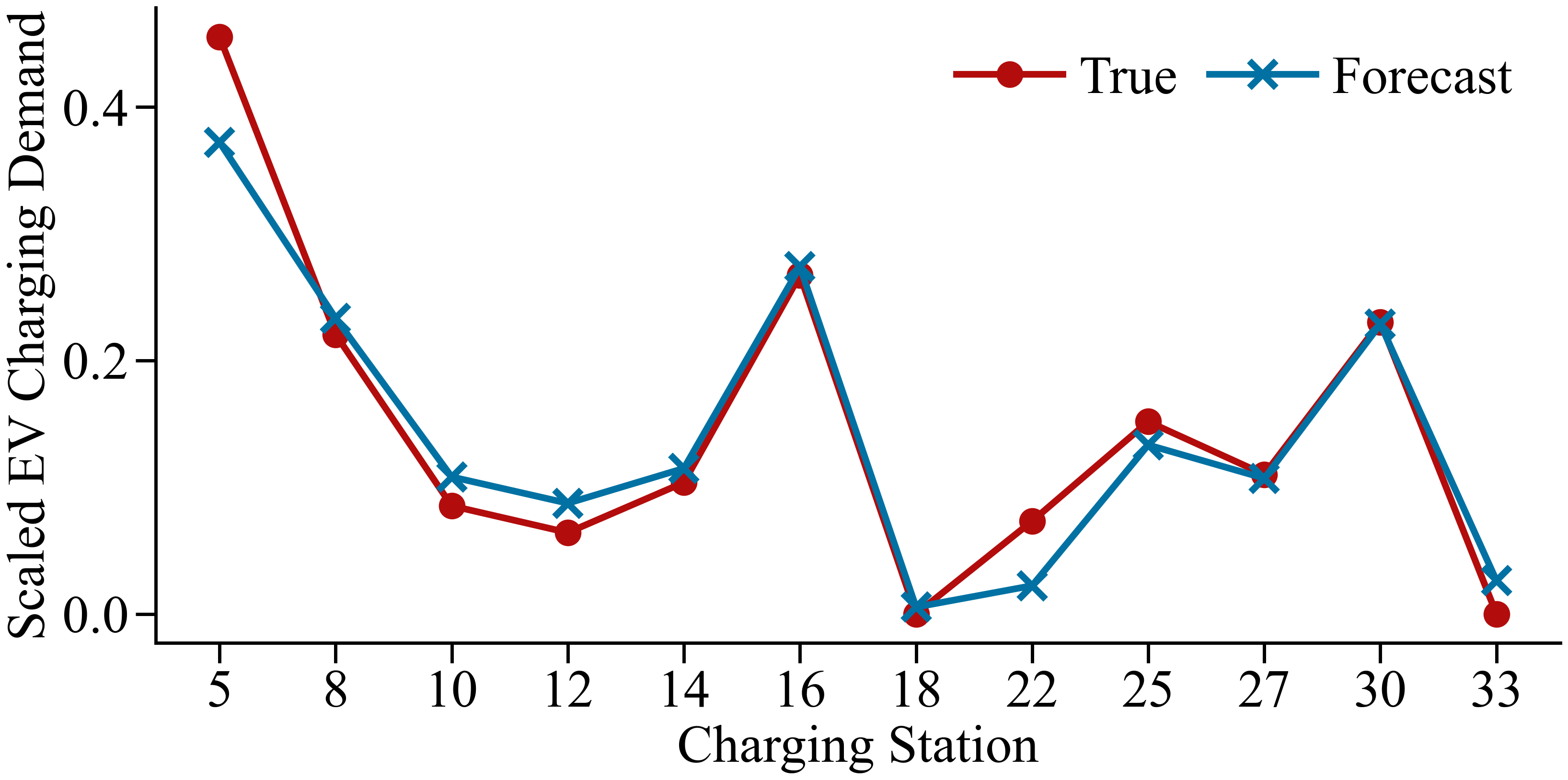}}
        \caption{}
        \label{ASTGCNpred1}
    \end{subfigure}
    \hfill
    \begin{subfigure}{\columnwidth}
        \setlength{\abovecaptionskip}{0.1cm}
        \setlength{\belowcaptionskip}{0.1cm}
        \centerline{\includegraphics[width=0.95\columnwidth]{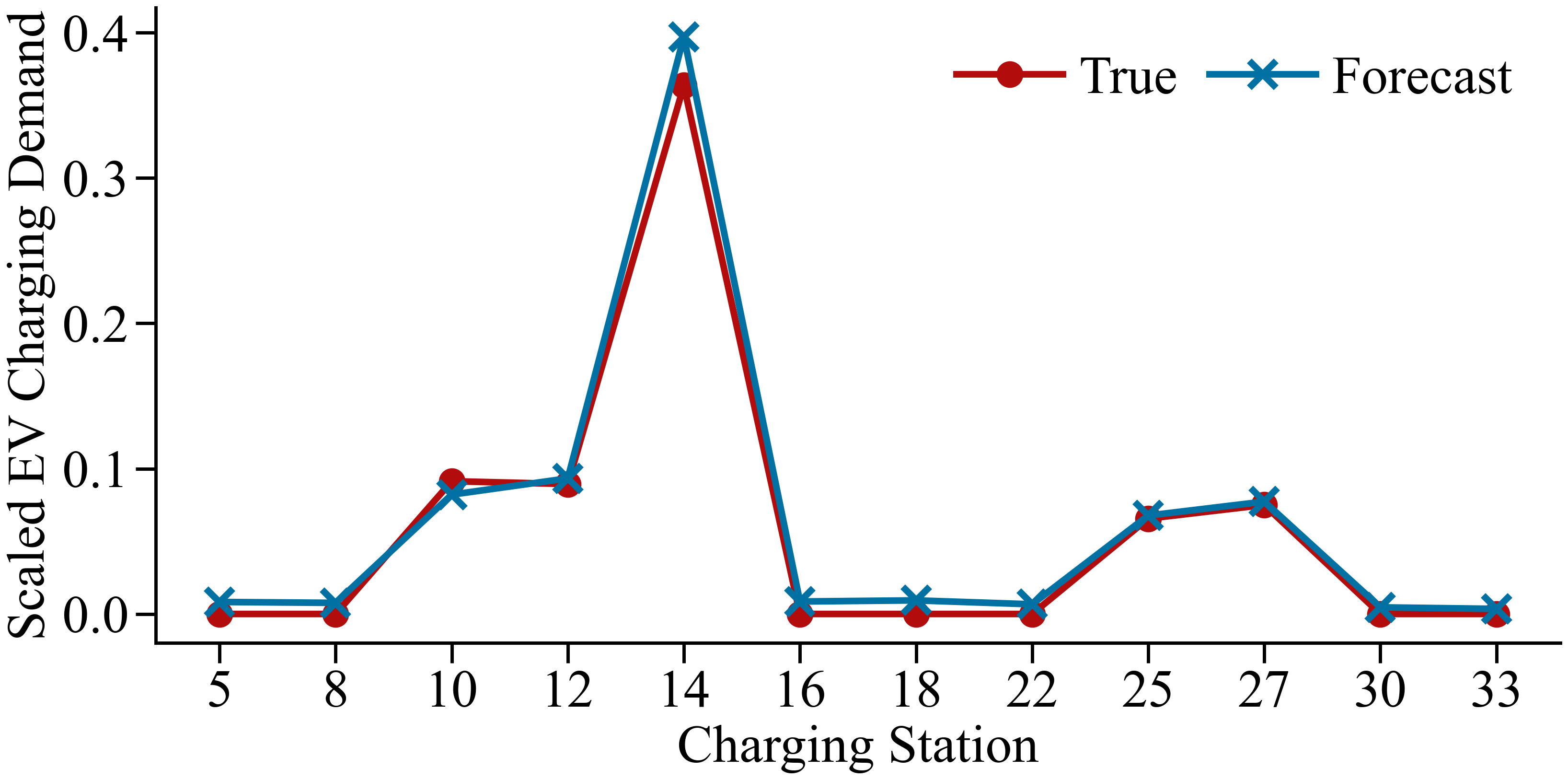}}
        \caption{}
        \label{ASTGCNpred2}
    \end{subfigure}
    \caption{\rmfamily Forecasting results of two test samples.}
    \label{ASTGCNpred}
\end{figure}

To visualize the deterministic forecasting results of ASTGCN,
we randomly select the forecasting results of two samples from the test set, 
as shown in Fig.~\ref{ASTGCNpred}.
The forecasted values (blue x marks) closely match the true values (red dots) with minimal forecasting bias, 
demonstrating the model's high forecasting precision. 
Furthermore, the peak charging demands observed at station 5 in Fig.~\ref{ASTGCNpred}.(a) and station 14 in Fig.~\ref{ASTGCNpred}.(b) 
are accurately forecasted. 
Similarly, the low charging demands at stations 18 and 33 in Fig.~\ref{ASTGCNpred}.(a) and multiple stations in Fig.~\ref{ASTGCNpred}.(b) 
are also well forecasted, 
indicating ASTGCN's robustness in handling both high and low charging demands.

\subsection{Probabilistic Forecasting Performance}
\begin{figure}[htbp]
    \setlength{\abovecaptionskip}{-0.1cm}
    \setlength{\belowcaptionskip}{-0.1cm}
    \centerline{\includegraphics[width=0.95\columnwidth]{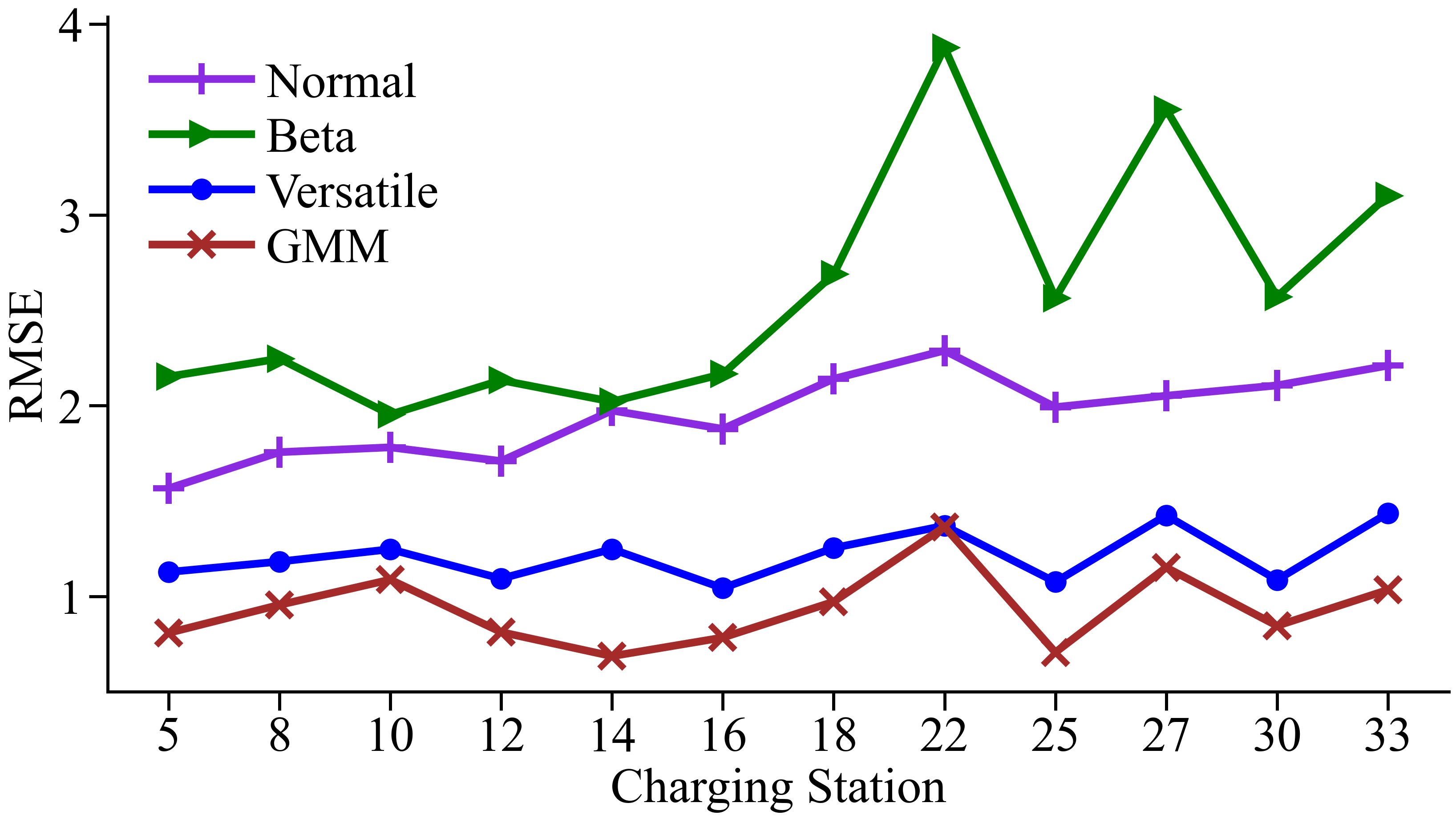}}
    \caption{\rmfamily RMSE comparison of different distributions.}
    \label{ErrorfitRMSE}
\end{figure}

Recall that we conduct probabilistic forecasting by distribution fitting (GMM in this paper) 
of forecasting errors based on predefined deterministic forecasting intervals.
We evaluate the performance of probabilistic forecasting by calculating the RMSE 
of the fitted distribution and the empirical distribution.
Besides, we also compare the performance of GMM with normal, beta, and versatile distributions~\cite{Versatile}.
Fig.~\ref{ErrorfitRMSE} presents the RMSE comparison of
the four distributions across the 12 charging stations.
Fig.~\ref{ErrorfitRMSE} reveals that GMM outperforms the other three distributions in terms of fitting accuracy across all stations. 
The consistently lower RMSE of GMM demonstrates its robustness and flexibility in fitting complex error distributions. 
While versatile distribution also shows improved performance, it still falls short compared to GMM. 
Normal and beta distributions, with their higher and more variable RMSE, 
are less effective in capturing the intricate data patterns. 
Therefore, to ensure accurate fitting of complex distributions, 
GMM is the preferred method.

\begin{figure}[htbp]
    \centering
    \begin{subfigure}{\columnwidth}
        \centerline{\includegraphics[width=0.8\columnwidth]{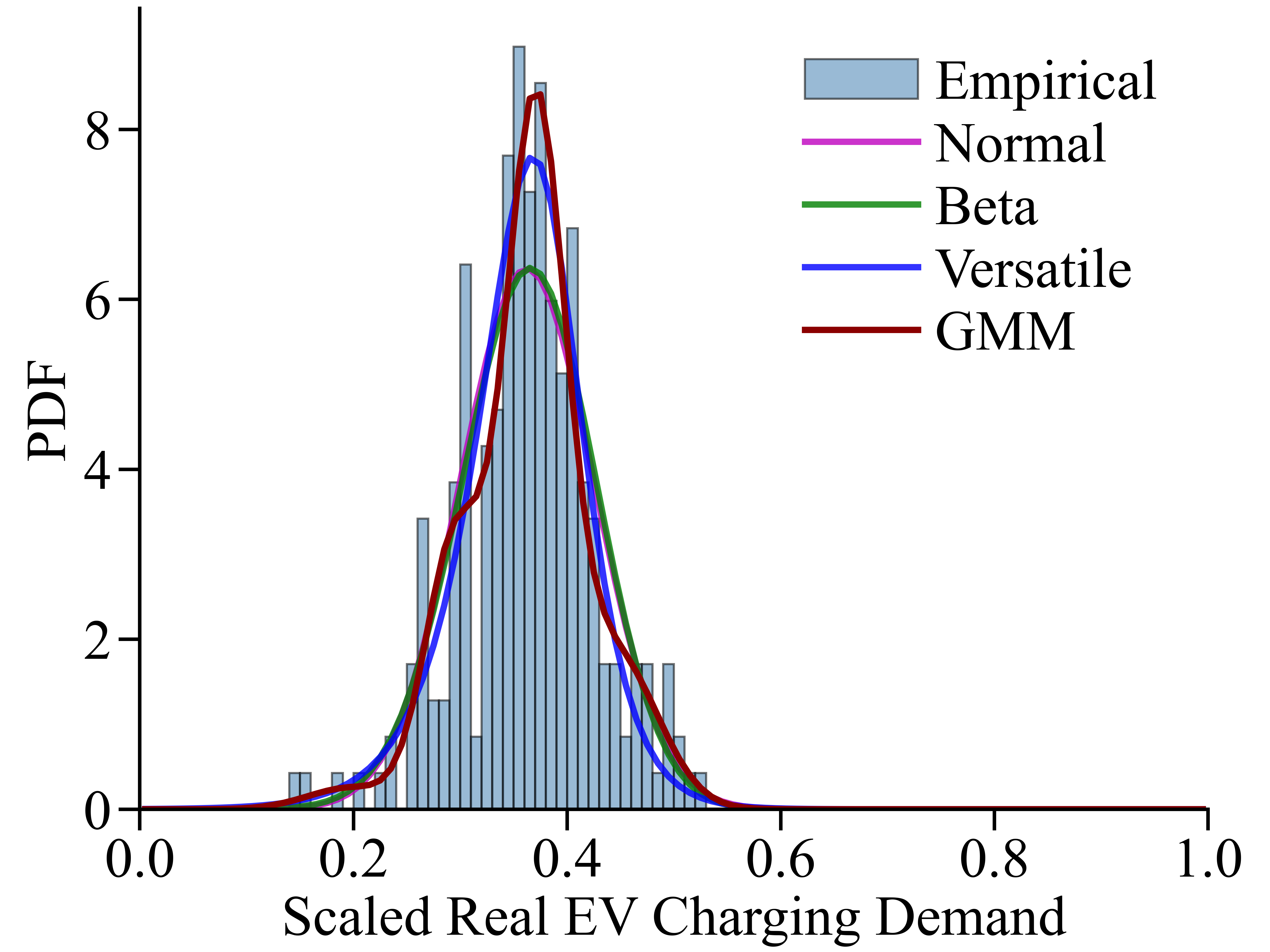}}
        \caption{}
        \label{Errorfitsub1}
    \end{subfigure}
    \hfill
    \begin{subfigure}{\columnwidth}
        \centerline{\includegraphics[width=0.8\columnwidth]{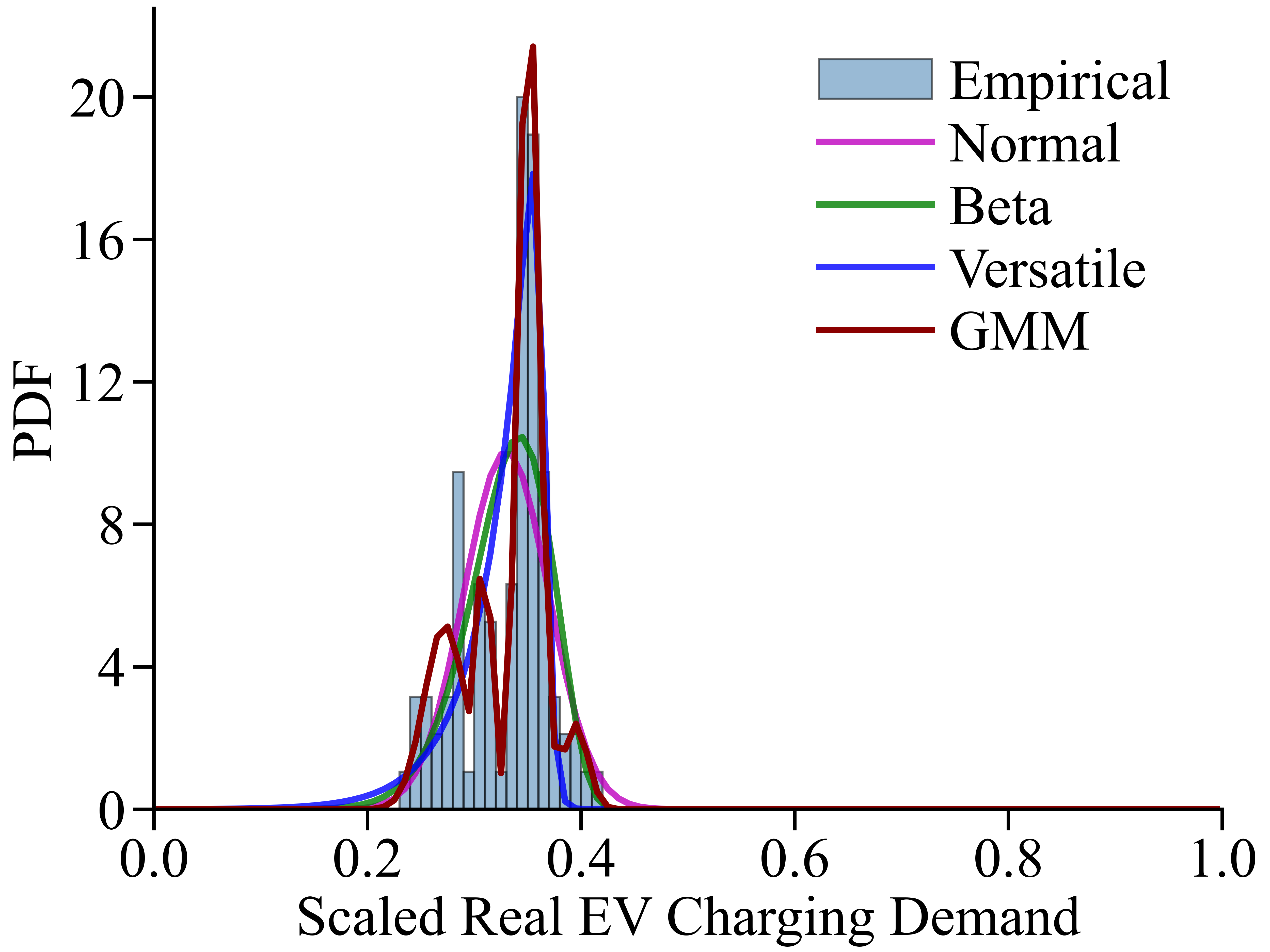}}
        \caption{}
        \label{Errorfitsub2}
    \end{subfigure}
    \hfill
    \begin{subfigure}{\columnwidth}
        \centerline{\includegraphics[width=0.8\columnwidth]{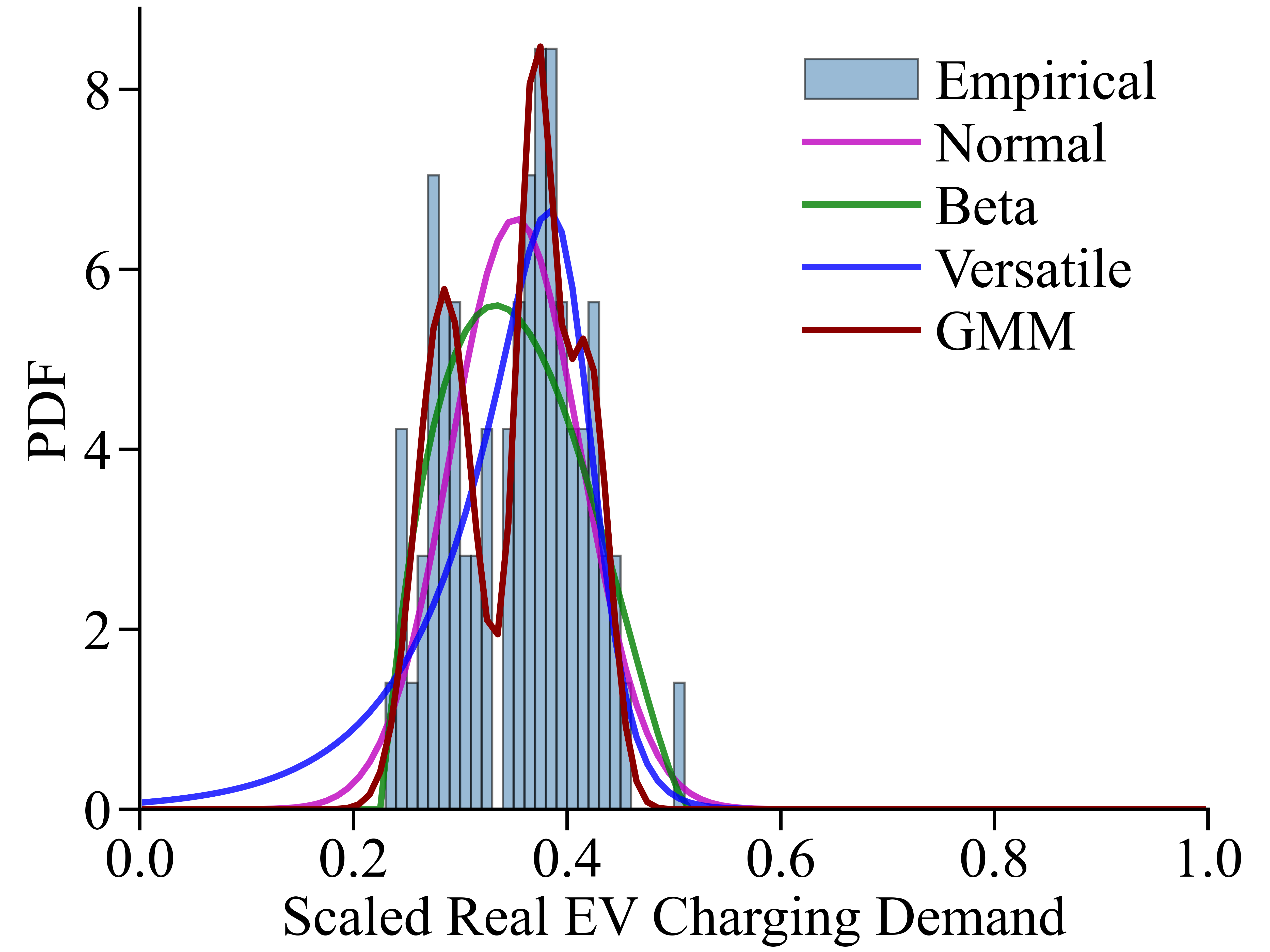}}
        \caption{}
        \label{Errorfitsub3}
    \end{subfigure}
    \caption{\rmfamily Fitting comparison of GMM, normal, beta and versatile: 
        (a) station 5, (b) station 14, (c) station 27.}
    \label{Errorfit}
\end{figure}

To visualize the fitting performance,
we select three stations (station 5, station 14, and station 27) and select the forecasting interval $[0.39,0.40]$ for comparison.
The fitting results of true values are shown in Fig.~\ref{Errorfit}.

The empirical distribution at station 5 in~Fig.~\ref{Errorfit}.(a) exhibits a clear unimodal and approximately symmetric shape. 
All four distributions demonstrate good fitting performance, accurately capturing the central peak. 
This suggests that for simple, symmetric distributions, all four distributions can provide a reasonable fit.
The empirical distribution at station 14 in~Fig.~\ref{Errorfit}.(b) is characterized by a very sharp unimodal peak. 
Here, versatile and GMM distributions show their adaptability to sharp features,
while normal and beta distributions fail to capture this sharpness, 
highlighting their limitations in representing distributions with steep gradients. 
The empirical distribution at station 27 in~Fig.~\ref{Errorfit}.(c) displays a complex multimodal distribution with multiple peaks. 
In this case, only the GMM is capable of accurately fitting the multiple peaks,
while the normal, beta, and versatile distributions expose their limitations in handling complex, multi-peaked data.

\subsection{Real-Time Risk Analysis Performance}

\begin{figure}[htbp]
    \setlength{\abovecaptionskip}{-0.1cm}
    \setlength{\belowcaptionskip}{-0.1cm}
    \centerline{\includegraphics[width=0.95\columnwidth]{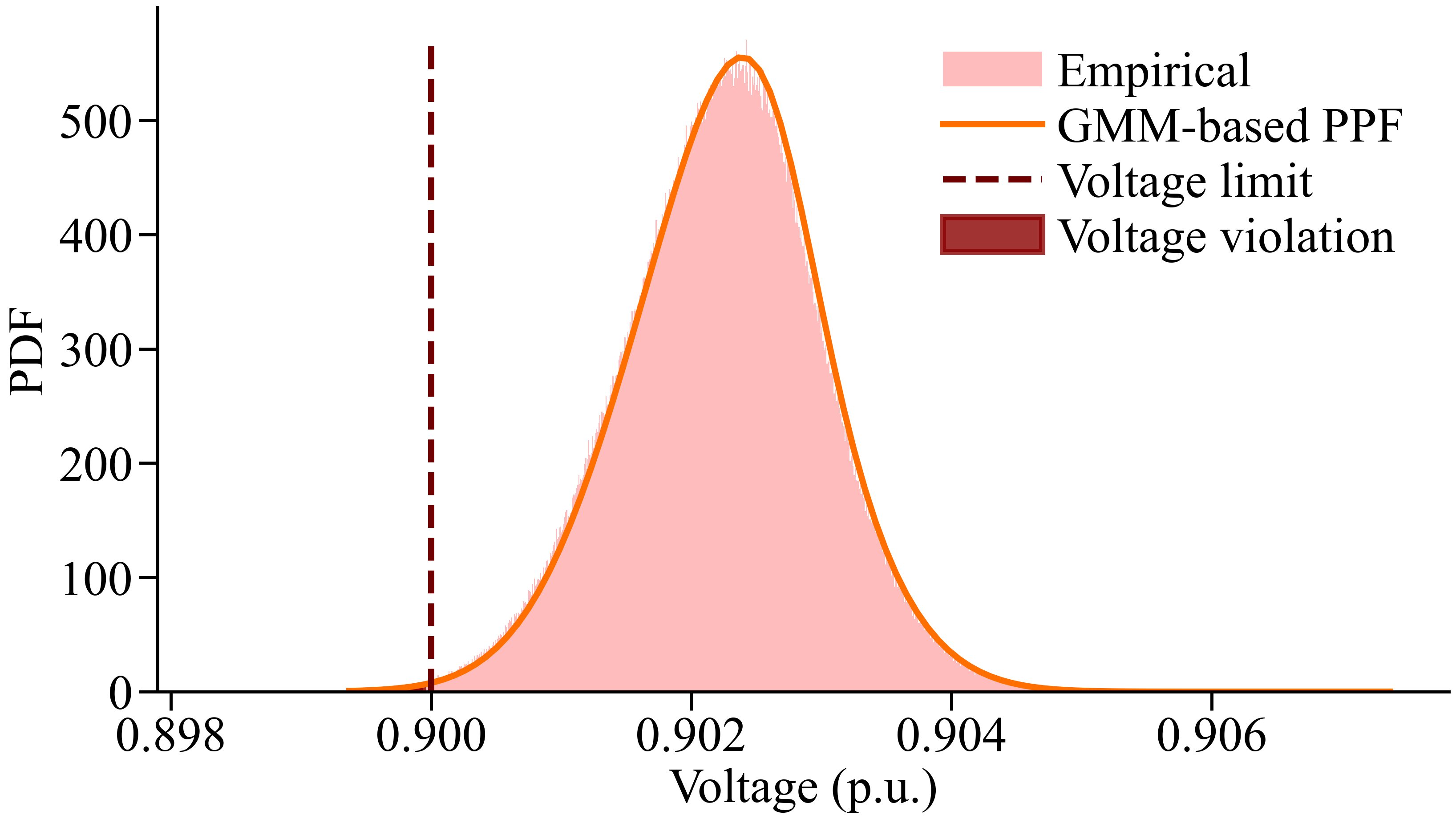}}
    \caption{\rmfamily Comparison of GMM-based PPF and MC-based PPF of bus 18.}
    \label{RiskAnalysis}
\end{figure}

In terms of PPF calculation,
we compare the performance of GMM-based PPF with the Monte Carlo simulation method (i.e., the benchmark).
Recall that we focus on voltage drop as the predominant risk factor.
Bus 18, located at the end of the IEEE 33 network,
faces the most significant voltage drop.
The PPF calculation result of bus 18 is illustrated in Fig.~\ref{RiskAnalysis}.
It is observed that GMM-based PPF (orange curve) closely aligns with the empirical probability density based on Monte Carlo simulation (pink area),
indicating that GMM-based PPF effectively captures the characteristics of the probabilistic voltage distribution.
Regarding the calculation time, Monte Carlo simulation takes over 15 minutes with $10^6$ samples,
while GMM-based PPF only requires around 4 seconds, achieving high efficiency in PPF calculation.
The dark red area represents the low voltage violation region. 
Without regulation, the system safety would be compromised, 
leading to severe adjustment costs,
highlighting the importance of risk analysis and HC assessment for EVs.

Regarding the operational boundary identification, 
we set $\varsigma =0.1\%$, and the corresponding voltage limit is $\underline{V}=0.9$ based on~\eqref{operational boundaries}.
The dark red area in Fig.~\ref{RiskAnalysis} denotes the region where the voltage drops below the operational boundary.

\subsection{Performance of Real-time HC assessment for EVs}
\begin{figure}[htbp]
    \setlength{\abovecaptionskip}{-0.1cm}  
    \setlength{\belowcaptionskip}{-0.1cm}   
    \centerline{\includegraphics[width=0.95\columnwidth]{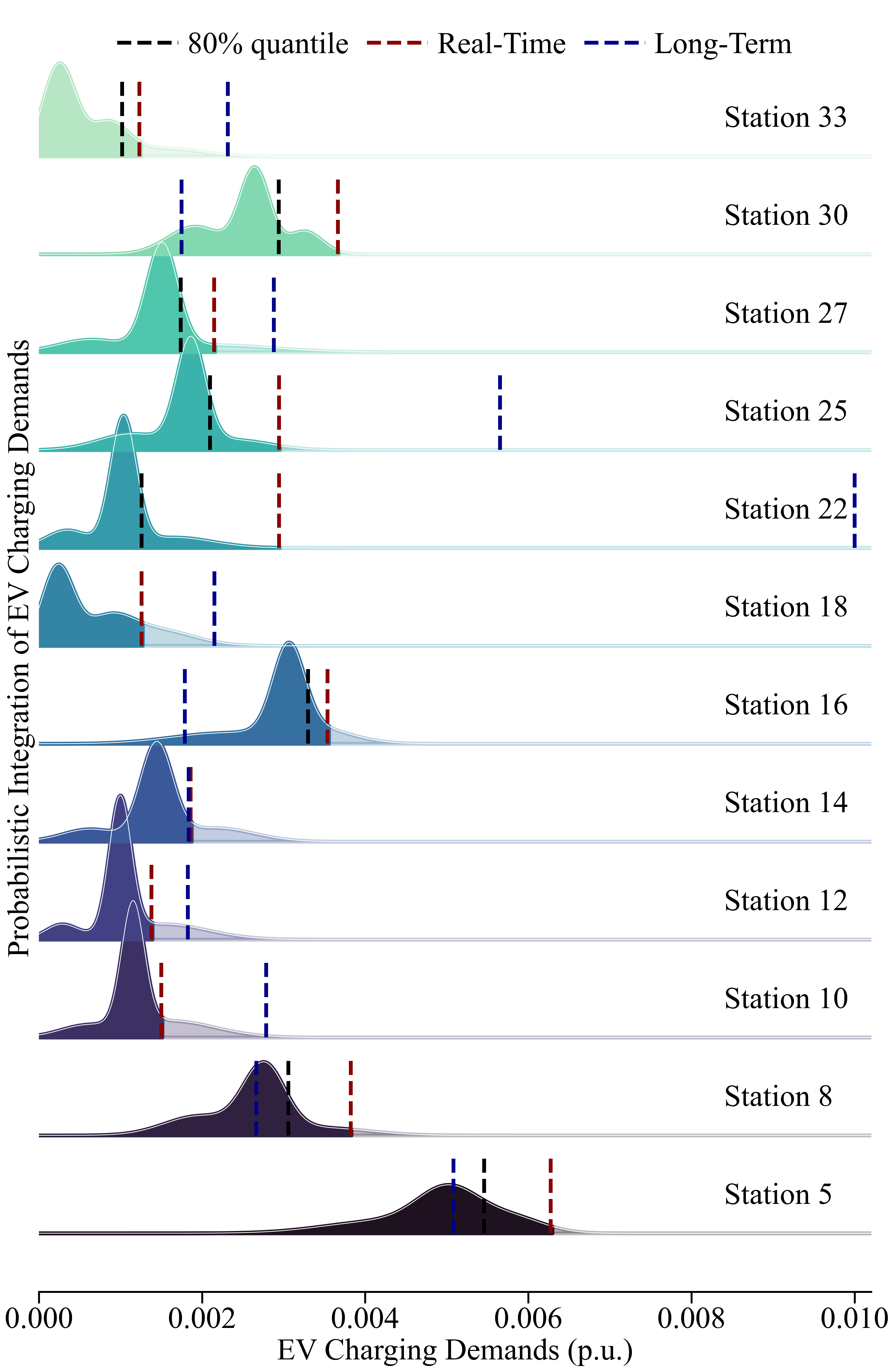}}
    \caption{\rmfamily Real-time and long-term HC of EVs.}
    \label{HCEV21injection}
\end{figure}

\begin{figure}[htbp]
    \setlength{\abovecaptionskip}{-0.1cm}  
    \setlength{\belowcaptionskip}{-0.1cm}   
    \centerline{\includegraphics[width=0.95\columnwidth]{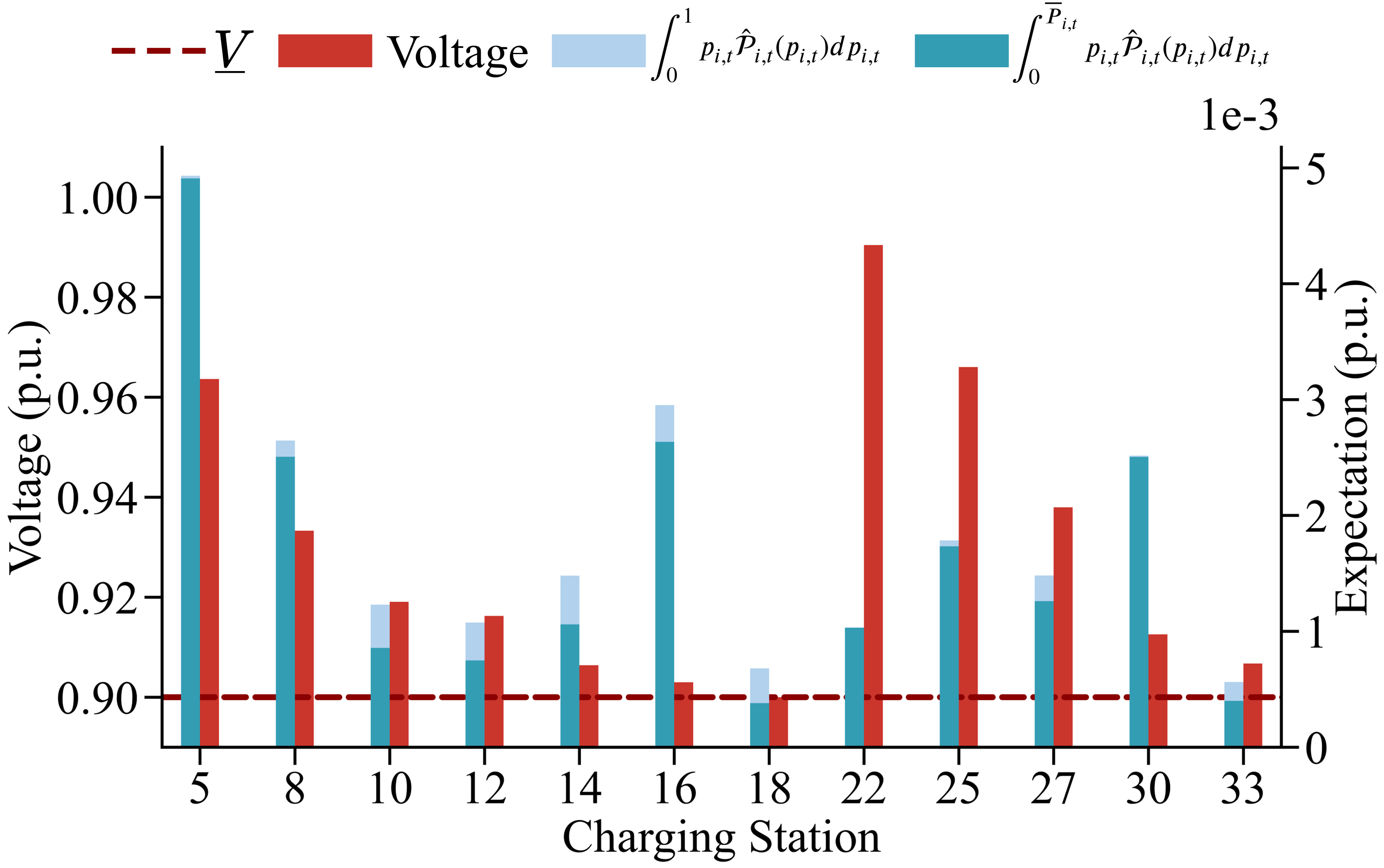}}
    \caption{\rmfamily Expected charging demands satisfaction and bus voltage of real-time HC assessment of EVs.}
    \label{HCEV21satisfaction}
\end{figure}

By solving the optimization problem in~\eqref{opt},
the real-time HC of EVs at each bus can be assessed in real-time ($<0.5$ seconds).
Setting $\{\epsilon_i\}_{i\in \varOmega^{\mathrm{C}} }=0.2$, 
the HC assessment results are illustrated in Fig.~\ref{HCEV21injection}
and Fig.~\ref{HCEV21satisfaction}.

In Fig.~\ref{HCEV21injection},
the filled area denotes the forecasted probabilistic distribution of EV charging demands,
the black dotted lines are the minimum satisfaction level of EV charging demands 
(i.e., $\hat{\mathcal{F}}_{i,t}^{-1}(1-\epsilon_i)$ for station $i$),
the red dotted lines denote $\overline{P}_{i,t}$ of real-time HC assessment,
while the blue dotted lines denote $\overline{P}_{i}$ of long-term HC assessment.
Notably, charging station 18 and charging station 33,
located at the ends of the distribution network, 
exhibit relatively low EV HC values. 
This observation is consistent with expectations, 
as end nodes typically experience more significant voltage drops due to their distal location in the network.
Conversely, charging stations 5, 22, 25 display relatively high EV HC values. 
This can be attributed to their greater voltage safety margins, 
allowing for greater flexibility in accommodating EV charging demands without violating voltage constraints. 
Additionally, the HC results from the long-term approach (blue dotted lines) 
show both overestimations (e.g., at charging stations 22 and 25) 
and underestimations (e.g., at charging stations 16 and 30), 
which fail to accurately align with the real-time probabilities of EV charging demands.
Moreover, the long-term estimated HC values tend to aggregate EV charging loads at a limited number of nodes 
with substantial voltage margins (i.e., stations 22 and 25), 
while excluding other nodes,
which is impractical in real-world traffic networks and EV system operations. 
The proposed method improves the real-time HC of EV to 0.2, 
a 64\% increase compared to the 0.122 of long-term assessment.

Fig.~\ref{HCEV21satisfaction} illustrates the voltage levels (red columns) at each station in $\varOmega^{\mathrm{C}}$,
the expected EV charging demands (light blue columns) 
and the expected satisfaction of EV charging demands (blue-green columns).
It is observed that the expected satisfaction of EV charging demands reaches a relatively low level at the buses 
with low voltage levels.
This is because the low voltage levels limit the EV charging demands that can be accommodated.
Additionally, the optimization problem~\eqref{opt} achieves an intricate balance to manage charging demands across the network,
and ensures that despite the variability in voltage levels, the network can sustain the charging demands to the greatest extent possible. 

\subsection{Performance Comparison of Deterministic Forecasting Models}
To validate the superiority of the proposed ASTGCN model, 
we compare the deterministic forecasting performance of ASTGCN with two categories of forecasting models, 
including statistical and deep learning models.
\begin{enumerate}
    \item Historical Average (HA). 
        HA forecasts future values by averaging the historical observations. 
    \item  Auto-regressive Integrated Moving Average (ARIMA) \cite{ARIMA}. 
        ARIMA assumes that time series can be made stationary through differencing. 
        Future values are forecasted based on a linear combination of past observations, 
        errors, and differenced values to account for trends.
    \item Multilayer Perceptron (MLP): 
        MLP is a neural network with multiple fully connected layers, 
        used for time series forecasting by identifying patterns in sequential data.
    \item Bi-directional Long Short-Term Memory (Bi-LSTM). 
        LSTM is a type of recurrent neural network designed to capture long-term dependencies. 
        BiLSTM enhances LSTM by processing data in both forward and backward directions, 
        capturing richer context from past and future sequences.
    \item Temporal Convolutional Network (TCN). 
        TCN is a type of convolutional neural network designed to capture temporal dependencies in time series. 
        TCN uses dilated convolutions to increase the receptive field of the network, 
        allowing them to capture long-term dependencies in time series.
    \item Combined Convolutional Neural Network and Bi-directional Long Short-Term Memory Network (CNN-BiLSTM).
        CNN-BiLSTM is a hybrid model that combines the strengths of convolutional neural networks (CNN) and BiLSTM. 
        CNN can extract variable coupling features within short time windows of multi-dimensional variables 
        and reduce the dimensionality of the input multi-dimensional data. 
        BiLSTM layer extracts forward and backward dependencies from long time series.
    \item Spatio-Temporal Graph Convolutional Network (STG\newline CN).
    STGCN is a type of graph convolutional network designed to capture both spatial and temporal dependencies in time series. 
    We adopt the model from~\cite{STGCN}, 
    in which a fixed graph structure is adopted
    using the adjacency matrix constructed by the Euclidean distance-based similarity between historical charging observations.
\end{enumerate}

The comparison results of the above forecasting models
under MAE, WAPE, and RMSE metrics are presented in Table.~\ref{PerformanceComparison}.
It is observed that the proposed ASTGCN model outperforms all other models across all metrics.
The statistical models (HA and ARIMA) exhibit unsatisfactory performance. 
Specifically, HA exhibits significantly higher MAE and RMSE values compared to other models, and ARIMA shows a notably higher WAPE.
This can be attributed to the inherent limitations of statistical models in capturing the complex, 
nonlinear spatial and temporal patterns within the data.
Statistical models rely on linear assumptions without accounting for the intricate dependencies that often characterize real-world datasets, 
leading to suboptimal forecasting accuracy.
In contrast, deep learning models (MLP, TCN, BiLSTM, CNN-BiLSTM, STGCN and ASTGCN) demonstrate markedly better performance. 
First, the time series models (TCN and BiLSTM) outperform MLP, 
highlighting the critical importance of capturing temporal correlations, 
which simple feed forward neural networks like MLP are unable to do effectively.
Moreover, CNN-BiLSTM, which integrates CNN to extract spatial features, 
achieves better performance than the standalone time series models. 
STGCN further improves the forecasting performance by constructing a graph structure to capture spatial dependencies.
However, its use of a fixed graph structure to model these dependencies neglects the time-varying nature. 
The proposed ASTGCN model addresses this limitation by constructing a combined time-invariant and time-varying graph structure 
through the learning of an adaptive weighted adjacency matrix. 
As a result, 
ASTGCN achieves the best performance across all metrics, 
demonstrating its effectiveness and superiority over both statistical models and other deep learning models.

\begin{table}[htbp]
    \centering
    \fontfamily{ptm}\selectfont 
    \setlength{\abovecaptionskip}{-0.1cm}
    \setlength{\belowcaptionskip}{-0.1cm}
    \caption{\rmfamily Performance comparison of deterministic forecasting models.}
    \label{PerformanceComparison}
    \begin{tabular}{ccccc}
        \toprule
        \textbf{Type} & \textbf{Model} & \textbf{MAE} & \textbf{WAPE (\%)} & \textbf{RMSE} \\
        \midrule
        \multirow{2}{*}{\makecell{Statistical\\ Models}} & HA & 0.0739 & 62.57 & 0.0959 \\
        & ARIMA & 0.0647 & 110.73 & 0.0783 \\
        \midrule
        \multirow{6}{*}{\makecell{Deep \\Learning \\Models}} & MLP & 0.0491 & 42.05 & 0.0600 \\
        & TCN & 0.0397 & 34.02 & 0.0577 \\
        & BiLSTM & 0.0454 & 38.89 & 0.0564 \\
        & CNN-BiLSTM & 0.0442 & 37.90 & 0.0548 \\
        & STGCN & 0.0330 & 28.69 & 0.0486 \\
        & ASTGCN & 0.0285 & 24.45 & 0.0442 \\
        \bottomrule
    \end{tabular}
\end{table}

\begin{figure}[htbp]
    \setlength{\abovecaptionskip}{-0.1cm}
    \setlength{\belowcaptionskip}{-0.1cm}
    \centerline{\includegraphics[width=0.95\columnwidth]{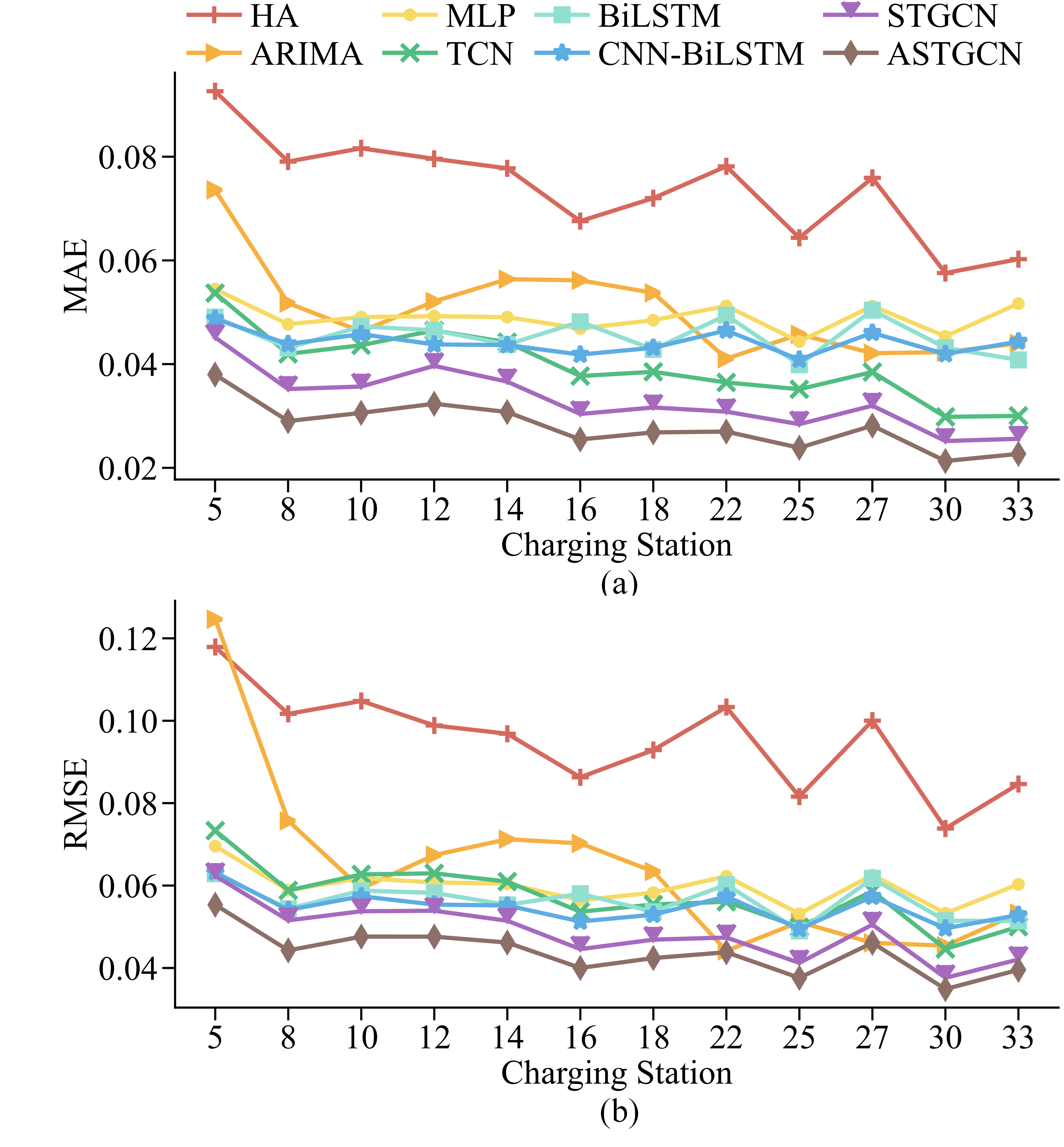}}
    \caption{\rmfamily Performance comparison of different forecasting models across 12 charging stations: (a) MAE, (b) RMSE.}
    \label{PerformanceComparisonnode}
\end{figure}

To further evaluate the performance of ASTGCN, 
we conduct a node-wise comparison, 
as depicted in Fig.~\ref{PerformanceComparisonnode}. 
The forecasting performance varies across charging stations, 
reflecting the distinct characteristics and complexities of the charging demands at each station. 
Despite these variations, 
a consistent performance tendency is observed across different charging stations.
Notably, ASTGCN consistently outperforms other models in both MAE and RMSE metrics.
Additionally, ASTGCN demonstrates remarkable stability in its performance across different charging stations. 
This stability is crucial for practical applications where consistent performance is required for reliable forecasting.
These results demonstrate that the proposed ASTGCN model is a powerful tool for charging demands forecasting 
of multiple charging stations.

\subsection{Ablation Study}
To investigate the effect of each module within the proposed ASTGCN model, 
we conduct an ablation study. 
Three variants of ASTGCN are as follows:
\begin{enumerate}
    \item ASTGCN-noWA:  This variant removes the adaptive graph module, 
            and adopts the similarity matrix of historical charging demands observations~\cite{STMGCN} 
            to describe the spatial correlation.
    \item ASTGCN-noTA: This variant removes the temporal attention module. 
    \item ASTGCN-fc:   This variant replaces the second-order pooling module with traditional 
                        fully connected layer.
\end{enumerate}

Table.~\ref{AblationStudy} presents the performance comparison of the three variants above with ASTGCN.
ASTGCN-noWA achieves the highest MAE (0.0338), WAPE (28.96\%) and RMSE (0.0477),
indicating that the adaptive graph module plays a crucial role in capturing the underlying spatial dependencies.
Compared with ASTGCN-noWA,
ASTGCN-noTA retains the adaptive graph mechanism but removes the temporal attention module,
and achieves lower MAE (0.0302), WAPE (25.85\%) and RMSE (0.0458).
This result indicates that the temporal attention mechanism improves the forecasting accuracy,
but the spatial adaptive mechanism has more significant impacts.
ASTGCN-fc achieves the lowest MAE (0.0290), WAPE (24.90\%) and RMSE (0.0449), 
which are slightly higher than ASTGCN.
This result demonstrates that the second-order pooling module is effective in enhancing graph representation
by extracting second order graph information.
The proposed ASTGCN model combines the adaptive graph mechanism, temporal attention mechanism, and second-order pooling module,
and achieves the best performance with the lowest MAE (0.0285), WAPE (24.45\%) and RMSE (0.0442).
These results underscore the combined effectiveness of these modules in enhancing the accuracy of EV charging demands deterministic forecasting.

\begin{table}[htbp]
    \rmfamily 
    \centering
    \caption{\rmfamily Performance comparison of ASTGCN variants.}\label{AblationStudy}
    \begin{tabular}{cccc}
        \toprule
        \textbf{Model} & \textbf{MAE} & \textbf{WAPE (\%)}  & \textbf{RMSE} \\
        \midrule
        ASTGCNnoWA & 0.0338 & 28.96 & 0.0477 \\
        ASTGCNnoTA & 0.0302 & 25.85 & 0.0458 \\
        ASTGCNfc & 0.0290 & 24.90 & 0.0449 \\
        ASTGCN & 0.0285 & 24.45 & 0.0442 \\
        \bottomrule
    \end{tabular}
\end{table}

\section{Conclusion}\label{conclusion}
In this paper, a real-time HC assessment method is proposed for EVs,
consisting of a three-step process of real-time EV charging demands probabilistic forecasting,
risk analysis and HC assessment.
In the first step,
the proposed ASTGCN is able to adaptively capture the spatial and temporal dependencies of multiple EV charging demands.
The comparison results of state-of-the-art forecasting models demonstrate the superiority of ASTGCN.
Meanwhile, the ablation study verifies the effectiveness of 
the adaptive graph mechanism in learning both time-invariant and time-varying graph patterns,
the temporal attention mechanism in capturing the relative temporal importance,
and the second-order pooling module in enhancing the graph representation.
In the second step,
with probabilistic forecasting the EV charging demands as GMM,
the proposed GMM-based PPF method is verified to describe the stochastic state of the power system and corresponding risk
with high efficiency.
In the third step, 
the proposed real-time HC formulation and assessment model are verified to effectively 
improve the HC of EVs compared to the long-term HC assessment.

In future work, we will explore improving the HC of EVs through smart coordination of
flexible energy resources (e.g., energy storage, demand response).
Also, by leveraging the flexibility of EVs, 
achieving smart charging strategies to enhance the overall sustainability of the power system
will be investigated.

\printcredits
\vspace{0.3cm}

\noindent \textbf{Declaration of competing interest}

The authors declare that they have no known competing financial interests or personal relationships that could have appeared to influence the work reported in this paper.
\vspace{0.3cm}

\noindent \textbf{Acknowledgements}

This work was supported in part by 
Special Foundation of Jiangsu Province Innovation Support Program (Soft Science Research) under Grant (No. BE2023093-1), National Natural Science Foundation of China under Grant (No. 52307144), the special funding from China Postdoctoral Science Foundation (No.2023TQ0169), and the Project (SKLD24KZ04) supported by State Key Laboratory of Power System Operation and Control.

\bibliographystyle{model1-num-names}

\bibliography{EVHC}

\begin{thebibliography}{65}
\expandafter\ifx\csname natexlab\endcsname\relax\def\natexlab#1{#1}\fi
\providecommand{\url}[1]{\texttt{#1}}
\providecommand{\href}[2]{#2}
\providecommand{\path}[1]{#1}
\providecommand{\DOIprefix}{doi:}
\providecommand{\ArXivprefix}{arXiv:}
\providecommand{\URLprefix}{URL: }
\providecommand{\Pubmedprefix}{pmid:}
\providecommand{\doi}[1]{\href{http://dx.doi.org/#1}{\path{#1}}}
\providecommand{\Pubmed}[1]{\href{pmid:#1}{\path{#1}}}
\providecommand{\bibinfo}[2]{#2}
\ifx\xfnm\relax \def\xfnm[#1]{\unskip,\space#1}\fi
\bibitem[{Lopes et~al.(2011)Lopes, Soares, and Almeida}]{EVimpact}
\bibinfo{author}{J.~A.~P. Lopes}, \bibinfo{author}{F.~J. Soares}, \bibinfo{author}{P.~M.~R. Almeida},
\newblock \bibinfo{title}{Integration of electric vehicles in the electric power system},
\newblock \bibinfo{journal}{Proceedings of the IEEE} \bibinfo{volume}{99} (\bibinfo{year}{2011}) \bibinfo{pages}{168--183}.
\bibitem[{IEA(2024)}]{iea2024EV}
\bibinfo{author}{IEA}, \bibinfo{title}{Global ev outlook 2024}, \bibinfo{howpublished}{\url{https://www.iea.org/reports/global-ev-outlook-2024}}, \bibinfo{year}{2024}.
\bibitem[{Clement-Nyns et~al.(2010)Clement-Nyns, Haesen, and Driesen}]{impact1}
\bibinfo{author}{K.~Clement-Nyns}, \bibinfo{author}{E.~Haesen}, \bibinfo{author}{J.~Driesen},
\newblock \bibinfo{title}{The impact of charging plug-in hybrid electric vehicles on a residential distribution grid},
\newblock \bibinfo{journal}{IEEE Transactions on Power Systems} \bibinfo{volume}{25} (\bibinfo{year}{2010}) \bibinfo{pages}{371--380}.
\bibitem[{Mu et~al.(2014)Mu, Wu, Jenkins, Jia, and Wang}]{OD}
\bibinfo{author}{Y.~Mu}, \bibinfo{author}{J.~Wu}, \bibinfo{author}{N.~Jenkins}, \bibinfo{author}{H.~Jia}, \bibinfo{author}{C.~Wang},
\newblock \bibinfo{title}{A spatial–temporal model for grid impact analysis of plug-in electric vehicles},
\newblock \bibinfo{journal}{Applied Energy} \bibinfo{volume}{114} (\bibinfo{year}{2014}) \bibinfo{pages}{456--465}.
\bibitem[{Xie et~al.(2021)Xie, Singh, Mitter, Dahleh, and Oren}]{xie2021toward}
\bibinfo{author}{L.~Xie}, \bibinfo{author}{C.~Singh}, \bibinfo{author}{S.~K. Mitter}, \bibinfo{author}{M.~A. Dahleh}, \bibinfo{author}{S.~S. Oren},
\newblock \bibinfo{title}{Toward carbon-neutral electricity and mobility: Is the grid infrastructure ready?},
\newblock \bibinfo{journal}{Joule} \bibinfo{volume}{5} (\bibinfo{year}{2021}) \bibinfo{pages}{1908--1913}.
\bibitem[{Ismael et~al.(2019)Ismael, {Abdel Aleem}, Abdelaziz, and Zobaa}]{HC}
\bibinfo{author}{S.~M. Ismael}, \bibinfo{author}{S.~H. {Abdel Aleem}}, \bibinfo{author}{A.~Y. Abdelaziz}, \bibinfo{author}{A.~F. Zobaa},
\newblock \bibinfo{title}{State-of-the-art of hosting capacity in modern power systems with distributed generation},
\newblock \bibinfo{journal}{Renewable Energy} \bibinfo{volume}{130} (\bibinfo{year}{2019}) \bibinfo{pages}{1002--1020}.
\bibitem[{Mulenga et~al.(2020)Mulenga, Bollen, and Etherden}]{HC_review}
\bibinfo{author}{E.~Mulenga}, \bibinfo{author}{M.~H. Bollen}, \bibinfo{author}{N.~Etherden},
\newblock \bibinfo{title}{A review of hosting capacity quantification methods for photovoltaics in low-voltage distribution grids},
\newblock \bibinfo{journal}{International Journal of Electrical Power \& Energy Systems} \bibinfo{volume}{115} (\bibinfo{year}{2020}) \bibinfo{pages}{105445}.
\bibitem[{Zenhom et~al.(2024)Zenhom, Aleem, Zobaa, and Boghdady}]{HC2}
\bibinfo{author}{Z.~M. Zenhom}, \bibinfo{author}{S.~H. E.~A. Aleem}, \bibinfo{author}{A.~F. Zobaa}, \bibinfo{author}{T.~A. Boghdady},
\newblock \bibinfo{title}{A comprehensive review of renewables and electric vehicles hosting capacity in active distribution networks},
\newblock \bibinfo{journal}{IEEE Access} \bibinfo{volume}{12} (\bibinfo{year}{2024}) \bibinfo{pages}{3672--3699}.
\bibitem[{Mousa et~al.(2024)Mousa, Mahmoud, and Lehtonen}]{HC3}
\bibinfo{author}{H.~H.~H. Mousa}, \bibinfo{author}{K.~Mahmoud}, \bibinfo{author}{M.~Lehtonen},
\newblock \bibinfo{title}{A comprehensive review on recent developments of hosting capacity estimation and optimization for active distribution networks},
\newblock \bibinfo{journal}{IEEE Access} \bibinfo{volume}{12} (\bibinfo{year}{2024}) \bibinfo{pages}{18545--18593}.
\bibitem[{Power(2019)}]{zero_corbon}
\bibinfo{author}{S.~Power}, \bibinfo{title}{Zero carbon communities}, \bibinfo{howpublished}{\url{https://www.spenergynetworks.co.uk/userfiles/file/Zero_Carbon_Communities_Report.pdf?v=3}}, \bibinfo{year}{2019}.
\bibitem[{Karmaker et~al.(2024)Karmaker, Prakash, Siddique, Hossain, and Pota}]{HCEVreview}
\bibinfo{author}{A.~K. Karmaker}, \bibinfo{author}{K.~Prakash}, \bibinfo{author}{M.~N.~I. Siddique}, \bibinfo{author}{M.~A. Hossain}, \bibinfo{author}{H.~Pota},
\newblock \bibinfo{title}{Electric vehicle hosting capacity analysis: Challenges and solutions},
\newblock \bibinfo{journal}{Renewable and Sustainable Energy Reviews} \bibinfo{volume}{189} (\bibinfo{year}{2024}) \bibinfo{pages}{113916}.
\bibitem[{Zhang et~al.(2020)Zhang, Yan, Liu, Zhang, and Lv}]{EVbehaviour}
\bibinfo{author}{J.~Zhang}, \bibinfo{author}{J.~Yan}, \bibinfo{author}{Y.~Liu}, \bibinfo{author}{H.~Zhang}, \bibinfo{author}{G.~Lv},
\newblock \bibinfo{title}{Daily electric vehicle charging load profiles considering demographics of vehicle users},
\newblock \bibinfo{journal}{Applied Energy} \bibinfo{volume}{274} (\bibinfo{year}{2020}) \bibinfo{pages}{115063}.
\bibitem[{{Pablo Carvallo} et~al.(2021){Pablo Carvallo}, Bieler, Collins, Mueller, Gehbauer, Gotham, and Larsen}]{EVimpact1}
\bibinfo{author}{J.~{Pablo Carvallo}}, \bibinfo{author}{S.~Bieler}, \bibinfo{author}{M.~Collins}, \bibinfo{author}{J.~Mueller}, \bibinfo{author}{C.~Gehbauer}, \bibinfo{author}{D.~J. Gotham}, \bibinfo{author}{P.~H. Larsen},
\newblock \bibinfo{title}{A framework to measure the technical, economic, and rate impacts of distributed solar, electric vehicles, and storage},
\newblock \bibinfo{journal}{Applied Energy} \bibinfo{volume}{297} (\bibinfo{year}{2021}) \bibinfo{pages}{117160}.
\bibitem[{Castelo~de Oliveira et~al.(2019)Castelo~de Oliveira, Bollen, Ribeiro, de~Carvalho, Zambroni, and Bonatto}]{DHC}
\bibinfo{author}{T.~E. Castelo~de Oliveira}, \bibinfo{author}{M.~Bollen}, \bibinfo{author}{P.~F. Ribeiro}, \bibinfo{author}{P.~M. de~Carvalho}, \bibinfo{author}{A.~C. Zambroni}, \bibinfo{author}{B.~D. J.~E. Bonatto},
\newblock \bibinfo{title}{The concept of dynamic hosting capacity for distributed energy resources: Analytics and practical considerations},
\newblock \bibinfo{journal}{Energies} \bibinfo{volume}{12} (\bibinfo{year}{2019}) \bibinfo{pages}{2576}.
\bibitem[{Zhang et~al.(2023)Zhang, Huang, Liao, and Liang}]{EV_tcn}
\bibinfo{author}{T.~Zhang}, \bibinfo{author}{Y.~Huang}, \bibinfo{author}{H.~Liao}, \bibinfo{author}{Y.~Liang},
\newblock \bibinfo{title}{A hybrid electric vehicle load classification and forecasting approach based on gbdt algorithm and temporal convolutional network},
\newblock \bibinfo{journal}{Applied Energy} \bibinfo{volume}{351} (\bibinfo{year}{2023}) \bibinfo{pages}{121768}.
\bibitem[{Zhang et~al.(2021)Zhang, Chan, Li, Wang, Qiu, and Wang}]{EVtf}
\bibinfo{author}{X.~Zhang}, \bibinfo{author}{K.~W. Chan}, \bibinfo{author}{H.~Li}, \bibinfo{author}{H.~Wang}, \bibinfo{author}{J.~Qiu}, \bibinfo{author}{G.~Wang},
\newblock \bibinfo{title}{Deep-learning-based probabilistic forecasting of electric vehicle charging load with a novel queuing model},
\newblock \bibinfo{journal}{IEEE Transactions on Cybernetics} \bibinfo{volume}{51} (\bibinfo{year}{2021}) \bibinfo{pages}{3157--3170}.
\bibitem[{Buzna et~al.(2021)Buzna, {De Falco}, Ferruzzi, Khormali, Proto, Refa, Straka, and {van der Poel}}]{prob_pred_EV2}
\bibinfo{author}{L.~Buzna}, \bibinfo{author}{P.~{De Falco}}, \bibinfo{author}{G.~Ferruzzi}, \bibinfo{author}{S.~Khormali}, \bibinfo{author}{D.~Proto}, \bibinfo{author}{N.~Refa}, \bibinfo{author}{M.~Straka}, \bibinfo{author}{G.~{van der Poel}},
\newblock \bibinfo{title}{An ensemble methodology for hierarchical probabilistic electric vehicle load forecasting at regular charging stations},
\newblock \bibinfo{journal}{Applied Energy} \bibinfo{volume}{283} (\bibinfo{year}{2021}) \bibinfo{pages}{116337}.
\bibitem[{Zhao et~al.(2020)Zhao, Wan, Song, and Cao}]{interval}
\bibinfo{author}{C.~Zhao}, \bibinfo{author}{C.~Wan}, \bibinfo{author}{Y.~Song}, \bibinfo{author}{Z.~Cao},
\newblock \bibinfo{title}{Optimal nonparametric prediction intervals of electricity load},
\newblock \bibinfo{journal}{IEEE Transactions on Power Systems} \bibinfo{volume}{35} (\bibinfo{year}{2020}) \bibinfo{pages}{2467--2470}.
\bibitem[{Li et~al.(2022)Li, Xu, Chew, Ding, and Zhao}]{prob_pred_PV}
\bibinfo{author}{Q.~Li}, \bibinfo{author}{Y.~Xu}, \bibinfo{author}{B.~S.~H. Chew}, \bibinfo{author}{H.~Ding}, \bibinfo{author}{G.~Zhao},
\newblock \bibinfo{title}{An integrated missing-data tolerant model for probabilistic pv power generation forecasting},
\newblock \bibinfo{journal}{IEEE Transactions on Power Systems} \bibinfo{volume}{37} (\bibinfo{year}{2022}) \bibinfo{pages}{4447--4459}.
\bibitem[{Khorramdel et~al.(2018)Khorramdel, Chung, Safari, and Price}]{kernal}
\bibinfo{author}{B.~Khorramdel}, \bibinfo{author}{C.~Y. Chung}, \bibinfo{author}{N.~Safari}, \bibinfo{author}{G.~C.~D. Price},
\newblock \bibinfo{title}{A fuzzy adaptive probabilistic wind power prediction framework using diffusion kernel density estimators},
\newblock \bibinfo{journal}{IEEE Transactions on Power Systems} \bibinfo{volume}{33} (\bibinfo{year}{2018}) \bibinfo{pages}{7109--7121}.
\bibitem[{Wan et~al.(2021)Wan, Cao, Lee, Song, and Ju}]{non_para_pred}
\bibinfo{author}{C.~Wan}, \bibinfo{author}{Z.~Cao}, \bibinfo{author}{W.-J. Lee}, \bibinfo{author}{Y.~Song}, \bibinfo{author}{P.~Ju},
\newblock \bibinfo{title}{An adaptive ensemble data driven approach for nonparametric probabilistic forecasting of electricity load},
\newblock \bibinfo{journal}{IEEE Transactions on Smart Grid} \bibinfo{volume}{12} (\bibinfo{year}{2021}) \bibinfo{pages}{5396--5408}.
\bibitem[{Zhang et~al.(2013)Zhang, Sun, Gao, Lin, and Cheng}]{Versatile}
\bibinfo{author}{Z.-S. Zhang}, \bibinfo{author}{Y.-Z. Sun}, \bibinfo{author}{D.~W. Gao}, \bibinfo{author}{J.~Lin}, \bibinfo{author}{L.~Cheng},
\newblock \bibinfo{title}{A versatile probability distribution model for wind power forecast errors and its application in economic dispatch},
\newblock \bibinfo{journal}{IEEE Transactions on Power Systems} \bibinfo{volume}{28} (\bibinfo{year}{2013}) \bibinfo{pages}{3114--3125}.
\bibitem[{Wan et~al.(2014)Wan, Xu, Pinson, Dong, and Wong}]{point_pred_not_reliable}
\bibinfo{author}{C.~Wan}, \bibinfo{author}{Z.~Xu}, \bibinfo{author}{P.~Pinson}, \bibinfo{author}{Z.~Y. Dong}, \bibinfo{author}{K.~P. Wong},
\newblock \bibinfo{title}{Probabilistic forecasting of wind power generation using extreme learning machine},
\newblock \bibinfo{journal}{IEEE Transactions on Power Systems} \bibinfo{volume}{29} (\bibinfo{year}{2014}) \bibinfo{pages}{1033--1044}.
\bibitem[{Zhang et~al.(2024)Zhang, Huang, Wang, Li, and Sun}]{tripchain}
\bibinfo{author}{L.~Zhang}, \bibinfo{author}{Z.~Huang}, \bibinfo{author}{Z.~Wang}, \bibinfo{author}{X.~Li}, \bibinfo{author}{F.~Sun},
\newblock \bibinfo{title}{An urban charging load forecasting model based on trip chain model for private passenger electric vehicles: A case study in beijing},
\newblock \bibinfo{journal}{Energy} \bibinfo{volume}{299} (\bibinfo{year}{2024}) \bibinfo{pages}{130844}.
\bibitem[{Dai et~al.(2014)Dai, Cai, Duan, and Zhao}]{EV_characteristic}
\bibinfo{author}{Q.~Dai}, \bibinfo{author}{T.~Cai}, \bibinfo{author}{S.~Duan}, \bibinfo{author}{F.~Zhao},
\newblock \bibinfo{title}{Stochastic modeling and forecasting of load demand for electric bus battery-swap station},
\newblock \bibinfo{journal}{IEEE Transactions on Power Delivery} \bibinfo{volume}{29} (\bibinfo{year}{2014}) \bibinfo{pages}{1909--1917}.
\bibitem[{Yan et~al.(2020)Yan, Zhang, Liu, Lv, Han, and Alfonzo}]{weather}
\bibinfo{author}{J.~Yan}, \bibinfo{author}{J.~Zhang}, \bibinfo{author}{Y.~Liu}, \bibinfo{author}{G.~Lv}, \bibinfo{author}{S.~Han}, \bibinfo{author}{I.~E.~G. Alfonzo},
\newblock \bibinfo{title}{Ev charging load simulation and forecasting considering traffic jam and weather to support the integration of renewables and evs},
\newblock \bibinfo{journal}{Renewable Energy} \bibinfo{volume}{159} (\bibinfo{year}{2020}) \bibinfo{pages}{623--641}.
\bibitem[{Pareschi et~al.(2020)Pareschi, Küng, Georges, and Boulouchos}]{survey}
\bibinfo{author}{G.~Pareschi}, \bibinfo{author}{L.~Küng}, \bibinfo{author}{G.~Georges}, \bibinfo{author}{K.~Boulouchos},
\newblock \bibinfo{title}{Are travel surveys a good basis for ev models? validation of simulated charging profiles against empirical data},
\newblock \bibinfo{journal}{Applied Energy} \bibinfo{volume}{275} (\bibinfo{year}{2020}) \bibinfo{pages}{115318}.
\bibitem[{Qian et~al.(2011)Qian, Zhou, Allan, and Yuan}]{MC1}
\bibinfo{author}{K.~Qian}, \bibinfo{author}{C.~Zhou}, \bibinfo{author}{M.~Allan}, \bibinfo{author}{Y.~Yuan},
\newblock \bibinfo{title}{Modeling of load demand due to ev battery charging in distribution systems},
\newblock \bibinfo{journal}{IEEE Transactions on Power Systems} \bibinfo{volume}{26} (\bibinfo{year}{2011}) \bibinfo{pages}{802--810}.
\bibitem[{Xiao et~al.(2020)Xiao, An, Cai, Wang, and Cai}]{EVquene}
\bibinfo{author}{D.~Xiao}, \bibinfo{author}{S.~An}, \bibinfo{author}{H.~Cai}, \bibinfo{author}{J.~Wang}, \bibinfo{author}{H.~Cai},
\newblock \bibinfo{title}{An optimization model for electric vehicle charging infrastructure planning considering queuing behavior with finite queue length},
\newblock \bibinfo{journal}{Journal of Energy Storage} \bibinfo{volume}{29} (\bibinfo{year}{2020}) \bibinfo{pages}{101317}.
\bibitem[{Arias et~al.(2017)Arias, Kim, and Bae}]{markovchain}
\bibinfo{author}{M.~B. Arias}, \bibinfo{author}{M.~Kim}, \bibinfo{author}{S.~Bae},
\newblock \bibinfo{title}{Prediction of electric vehicle charging-power demand in realistic urban traffic networks},
\newblock \bibinfo{journal}{Applied Energy} \bibinfo{volume}{195} (\bibinfo{year}{2017}) \bibinfo{pages}{738--753}.
\bibitem[{Amini et~al.(2016)Amini, Kargarian, and Karabasoglu}]{ARIMA}
\bibinfo{author}{M.~H. Amini}, \bibinfo{author}{A.~Kargarian}, \bibinfo{author}{O.~Karabasoglu},
\newblock \bibinfo{title}{Arima-based decoupled time series forecasting of electric vehicle charging demand for stochastic power system operation},
\newblock \bibinfo{journal}{Electric Power Systems Research} \bibinfo{volume}{140} (\bibinfo{year}{2016}) \bibinfo{pages}{378--390}.
\bibitem[{Wu et~al.(2022)Wu, Gu, Meng, Wen, and Ma}]{PCMP2}
\bibinfo{author}{K.~Wu}, \bibinfo{author}{J.~Gu}, \bibinfo{author}{L.~Meng}, \bibinfo{author}{H.~Wen}, \bibinfo{author}{J.~Ma},
\newblock \bibinfo{title}{An explainable framework for load forecasting of a regional integrated energy system based on coupled features and multi-task learning},
\newblock \bibinfo{journal}{Protection and Control of Modern Power Systems} \bibinfo{volume}{7} (\bibinfo{year}{2022}) \bibinfo{pages}{1--14}.
\bibitem[{Majidpour et~al.(2016)Majidpour, Qiu, Chu, Pota, and Gadh}]{EV_pred1}
\bibinfo{author}{M.~Majidpour}, \bibinfo{author}{C.~Qiu}, \bibinfo{author}{P.~Chu}, \bibinfo{author}{H.~R. Pota}, \bibinfo{author}{R.~Gadh},
\newblock \bibinfo{title}{Forecasting the ev charging load based on customer profile or station measurement?},
\newblock \bibinfo{journal}{Applied Energy} \bibinfo{volume}{163} (\bibinfo{year}{2016}) \bibinfo{pages}{134--141}.
\bibitem[{Wang et~al.(2023)Wang, Zhuge, Shao, Wang, Yang, and Wang}]{EV_lstm}
\bibinfo{author}{S.~Wang}, \bibinfo{author}{C.~Zhuge}, \bibinfo{author}{C.~Shao}, \bibinfo{author}{P.~Wang}, \bibinfo{author}{X.~Yang}, \bibinfo{author}{S.~Wang},
\newblock \bibinfo{title}{Short-term electric vehicle charging demand prediction: A deep learning approach},
\newblock \bibinfo{journal}{Applied Energy} \bibinfo{volume}{340} (\bibinfo{year}{2023}) \bibinfo{pages}{121032}.
\bibitem[{Bampos et~al.(2024)Bampos, Laitsos, Afentoulis, Vagropoulos, and Biskas}]{EV_pred3}
\bibinfo{author}{Z.~N. Bampos}, \bibinfo{author}{V.~M. Laitsos}, \bibinfo{author}{K.~D. Afentoulis}, \bibinfo{author}{S.~I. Vagropoulos}, \bibinfo{author}{P.~N. Biskas},
\newblock \bibinfo{title}{Electric vehicles load forecasting for day-ahead market participation using machine and deep learning methods},
\newblock \bibinfo{journal}{Applied Energy} \bibinfo{volume}{360} (\bibinfo{year}{2024}) \bibinfo{pages}{122801}.
\bibitem[{Shi et~al.(2024)Shi, Zhang, Bao, Gao, and Wang}]{STMGCN}
\bibinfo{author}{J.~Shi}, \bibinfo{author}{W.~Zhang}, \bibinfo{author}{Y.~Bao}, \bibinfo{author}{D.~W. Gao}, \bibinfo{author}{Z.~Wang},
\newblock \bibinfo{title}{Load forecasting of electric vehicle charging stations: Attention based spatiotemporal multi–graph convolutional networks},
\newblock \bibinfo{journal}{IEEE Transactions on Smart Grid} \bibinfo{volume}{15} (\bibinfo{year}{2024}) \bibinfo{pages}{3016--3027}.
\bibitem[{Ramadhani et~al.(2020)Ramadhani, Shepero, Munkhammar, Widén, and Etherden}]{PPFreview}
\bibinfo{author}{U.~H. Ramadhani}, \bibinfo{author}{M.~Shepero}, \bibinfo{author}{J.~Munkhammar}, \bibinfo{author}{J.~Widén}, \bibinfo{author}{N.~Etherden},
\newblock \bibinfo{title}{Review of probabilistic load flow approaches for power distribution systems with photovoltaic generation and electric vehicle charging},
\newblock \bibinfo{journal}{International Journal of Electrical Power \& Energy Systems} \bibinfo{volume}{120} (\bibinfo{year}{2020}) \bibinfo{pages}{106003}.
\bibitem[{Torquato et~al.(2018)Torquato, Salles, Pereira, Meira, and Freitas}]{PVHC1}
\bibinfo{author}{R.~Torquato}, \bibinfo{author}{D.~Salles}, \bibinfo{author}{C.~O. Pereira}, \bibinfo{author}{P.~C.~M. Meira}, \bibinfo{author}{W.~Freitas},
\newblock \bibinfo{title}{A comprehensive assessment of pv hosting capacity on low-voltage distribution systems},
\newblock \bibinfo{journal}{IEEE Transactions on Power Delivery} \bibinfo{volume}{33} (\bibinfo{year}{2018}) \bibinfo{pages}{1002--1012}.
\bibitem[{Ding and Mather(2016)}]{PVHC2}
\bibinfo{author}{F.~Ding}, \bibinfo{author}{B.~Mather},
\newblock \bibinfo{title}{On distributed pv hosting capacity estimation, sensitivity study, and improvement},
\newblock \bibinfo{journal}{IEEE Transactions on Sustainable Energy} \bibinfo{volume}{8} (\bibinfo{year}{2016}) \bibinfo{pages}{1010--1020}.
\bibitem[{Nijhuis et~al.(2017)Nijhuis, Gibescu, and Cobben}]{PPFGMM}
\bibinfo{author}{M.~Nijhuis}, \bibinfo{author}{M.~Gibescu}, \bibinfo{author}{S.~Cobben},
\newblock \bibinfo{title}{Gaussian mixture based probabilistic load flow for lv-network planning},
\newblock \bibinfo{journal}{IEEE Transactions on Power Systems} \bibinfo{volume}{32} (\bibinfo{year}{2017}) \bibinfo{pages}{2878--2886}.
\bibitem[{Barbosa et~al.(2020)Barbosa, Andrade, Torquato, Freitas, and Trindade}]{HCMC}
\bibinfo{author}{T.~Barbosa}, \bibinfo{author}{J.~Andrade}, \bibinfo{author}{R.~Torquato}, \bibinfo{author}{W.~Freitas}, \bibinfo{author}{F.~C. Trindade},
\newblock \bibinfo{title}{Use of ev hosting capacity for management of low-voltage distribution systems},
\newblock \bibinfo{journal}{IET Generation, Transmission \& Distribution} \bibinfo{volume}{14} (\bibinfo{year}{2020}) \bibinfo{pages}{2620--2629}.
\bibitem[{Baker(2019)}]{pf_fast1}
\bibinfo{author}{K.~Baker},
\newblock \bibinfo{title}{Learning warm-start points for ac optimal power flow},
\newblock in: \bibinfo{booktitle}{2019 IEEE 29th International Workshop on Machine Learning for Signal Processing (MLSP)}, \bibinfo{year}{2019}, pp. \bibinfo{pages}{1--6}.
\bibitem[{Donon et~al.(2020)Donon, Clément, Donnot, Marot, Guyon, and Schoenauer}]{pf_fast4}
\bibinfo{author}{B.~Donon}, \bibinfo{author}{R.~Clément}, \bibinfo{author}{B.~Donnot}, \bibinfo{author}{A.~Marot}, \bibinfo{author}{I.~Guyon}, \bibinfo{author}{M.~Schoenauer},
\newblock \bibinfo{title}{Neural networks for power flow: Graph neural solver},
\newblock \bibinfo{journal}{Electric Power Systems Research} \bibinfo{volume}{189} (\bibinfo{year}{2020}) \bibinfo{pages}{106547}.
\bibitem[{Zhao et~al.(2017)Zhao, Wang, Xu, Wang, Wan, and Chen}]{EVHC1}
\bibinfo{author}{J.~Zhao}, \bibinfo{author}{J.~Wang}, \bibinfo{author}{Z.~Xu}, \bibinfo{author}{C.~Wang}, \bibinfo{author}{C.~Wan}, \bibinfo{author}{C.~Chen},
\newblock \bibinfo{title}{Distribution network electric vehicle hosting capacity maximization: A chargeable region optimization model},
\newblock \bibinfo{journal}{IEEE Transactions on Power Systems} \bibinfo{volume}{32} (\bibinfo{year}{2017}) \bibinfo{pages}{4119--4130}.
\bibitem[{Alturki et~al.(2018)Alturki, Khodaei, Paaso, and Bahramirad}]{HCoptimization}
\bibinfo{author}{M.~Alturki}, \bibinfo{author}{A.~Khodaei}, \bibinfo{author}{A.~Paaso}, \bibinfo{author}{S.~Bahramirad},
\newblock \bibinfo{title}{Optimization-based distribution grid hosting capacity calculations},
\newblock \bibinfo{journal}{Applied Energy} \bibinfo{volume}{219} (\bibinfo{year}{2018}) \bibinfo{pages}{350--360}.
\bibitem[{Taheri et~al.(2021)Taheri, Jalali, Kekatos, and Tong}]{HC_tonglang}
\bibinfo{author}{S.~Taheri}, \bibinfo{author}{M.~Jalali}, \bibinfo{author}{V.~Kekatos}, \bibinfo{author}{L.~Tong},
\newblock \bibinfo{title}{Fast probabilistic hosting capacity analysis for active distribution systems},
\newblock \bibinfo{journal}{IEEE Transactions on Smart Grid} \bibinfo{volume}{12} (\bibinfo{year}{2021}) \bibinfo{pages}{2000--2012}.
\bibitem[{de~Lima et~al.(2023)de~Lima, Soares, Lezama, Franco, and Vale}]{EVoperation3}
\bibinfo{author}{T.~D. de~Lima}, \bibinfo{author}{J.~Soares}, \bibinfo{author}{F.~Lezama}, \bibinfo{author}{J.~F. Franco}, \bibinfo{author}{Z.~Vale},
\newblock \bibinfo{title}{A risk-based planning approach for sustainable distribution systems considering ev charging stations and carbon taxes},
\newblock \bibinfo{journal}{IEEE Transactions on Sustainable Energy} \bibinfo{volume}{14} (\bibinfo{year}{2023}) \bibinfo{pages}{2294--2307}.
\bibitem[{Abbasi et~al.(2019)Abbasi, Taki, Rajabi, Li, and Zhang}]{EVoperation4}
\bibinfo{author}{M.~H. Abbasi}, \bibinfo{author}{M.~Taki}, \bibinfo{author}{A.~Rajabi}, \bibinfo{author}{L.~Li}, \bibinfo{author}{J.~Zhang},
\newblock \bibinfo{title}{Coordinated operation of electric vehicle charging and wind power generation as a virtual power plant: A multi-stage risk constrained approach},
\newblock \bibinfo{journal}{Applied Energy} \bibinfo{volume}{239} (\bibinfo{year}{2019}) \bibinfo{pages}{1294--1307}.
\bibitem[{Wu et~al.(2019)Wu, Pan, Long, Jiang, and Zhang}]{wavenet}
\bibinfo{author}{Z.~Wu}, \bibinfo{author}{S.~Pan}, \bibinfo{author}{G.~Long}, \bibinfo{author}{J.~Jiang}, \bibinfo{author}{C.~Zhang},
\newblock \bibinfo{title}{Graph wavenet for deep spatial-temporal graph modeling},
\newblock \bibinfo{journal}{arXiv preprint arXiv:1906.00121}  (\bibinfo{year}{2019}).
\bibitem[{Kipf and Welling(2016)}]{spectralGCN}
\bibinfo{author}{T.~N. Kipf}, \bibinfo{author}{M.~Welling},
\newblock \bibinfo{title}{Semi-supervised classification with graph convolutional networks},
\newblock \bibinfo{journal}{arXiv preprint arXiv:1609.02907}  (\bibinfo{year}{2016}).
\bibitem[{Wu et~al.(2020)Wu, Pan, Chen, Long, Zhang, and Philip}]{wu2020comprehensive}
\bibinfo{author}{Z.~Wu}, \bibinfo{author}{S.~Pan}, \bibinfo{author}{F.~Chen}, \bibinfo{author}{G.~Long}, \bibinfo{author}{C.~Zhang}, \bibinfo{author}{S.~Y. Philip},
\newblock \bibinfo{title}{A comprehensive survey on graph neural networks},
\newblock \bibinfo{journal}{IEEE transactions on neural networks and learning systems} \bibinfo{volume}{32} (\bibinfo{year}{2020}) \bibinfo{pages}{4--24}.
\bibitem[{Dauphin et~al.(2017)Dauphin, Fan, Auli, and Grangier}]{Gated_CNN}
\bibinfo{author}{Y.~N. Dauphin}, \bibinfo{author}{A.~Fan}, \bibinfo{author}{M.~Auli}, \bibinfo{author}{D.~Grangier},
\newblock \bibinfo{title}{Language modeling with gated convolutional networks},
\newblock in: \bibinfo{booktitle}{International conference on machine learning}, \bibinfo{organization}{PMLR}, \bibinfo{year}{2017}, pp. \bibinfo{pages}{933--941}.
\bibitem[{Wang and Ji(2023)}]{SOP}
\bibinfo{author}{Z.~Wang}, \bibinfo{author}{S.~Ji},
\newblock \bibinfo{title}{Second-order pooling for graph neural networks},
\newblock \bibinfo{journal}{IEEE Transactions on Pattern Analysis and Machine Intelligence} \bibinfo{volume}{45} (\bibinfo{year}{2023}) \bibinfo{pages}{6870--6880}.
\bibitem[{Zhao et~al.(2017)Zhao, Wan, Xu, and Wang}]{EVoperation1}
\bibinfo{author}{J.~Zhao}, \bibinfo{author}{C.~Wan}, \bibinfo{author}{Z.~Xu}, \bibinfo{author}{J.~Wang},
\newblock \bibinfo{title}{Risk-based day-ahead scheduling of electric vehicle aggregator using information gap decision theory},
\newblock \bibinfo{journal}{IEEE Transactions on Smart Grid} \bibinfo{volume}{8} (\bibinfo{year}{2017}) \bibinfo{pages}{1609--1618}.
\bibitem[{Giraldo et~al.(2023)Giraldo, Arias, Vergara, Vlasiou, Hoogsteen, and Hurink}]{EVoperation2}
\bibinfo{author}{J.~S. Giraldo}, \bibinfo{author}{N.~B. Arias}, \bibinfo{author}{P.~P. Vergara}, \bibinfo{author}{M.~Vlasiou}, \bibinfo{author}{G.~Hoogsteen}, \bibinfo{author}{J.~L. Hurink},
\newblock \bibinfo{title}{Estimating risk-aware flexibility areas for electric vehicle charging pools via ac stochastic optimal power flow},
\newblock \bibinfo{journal}{Journal of Modern Power Systems and Clean Energy} \bibinfo{volume}{11} (\bibinfo{year}{2023}) \bibinfo{pages}{1247--1256}.
\bibitem[{Yan et~al.(2024)Yan, Huang, and Chen}]{EVAggregate}
\bibinfo{author}{D.~Yan}, \bibinfo{author}{S.~Huang}, \bibinfo{author}{Y.~Chen},
\newblock \bibinfo{title}{Real-time feedback based online aggregate ev power flexibility characterization},
\newblock \bibinfo{journal}{IEEE Transactions on Sustainable Energy} \bibinfo{volume}{15} (\bibinfo{year}{2024}) \bibinfo{pages}{658--673}.
\bibitem[{Madavan et~al.(2023)Madavan, Dahlin, Bose, and Tong}]{HCTL}
\bibinfo{author}{A.~N. Madavan}, \bibinfo{author}{N.~Dahlin}, \bibinfo{author}{S.~Bose}, \bibinfo{author}{L.~Tong},
\newblock \bibinfo{title}{Risk-based hosting capacity analysis in distribution systems},
\newblock \bibinfo{journal}{IEEE Transactions on Power Systems}  (\bibinfo{year}{2023}).
\bibitem[{Lin et~al.(2020)Lin, Bie, Pan, and Liu}]{linear_pf1}
\bibinfo{author}{C.~Lin}, \bibinfo{author}{Z.~Bie}, \bibinfo{author}{C.~Pan}, \bibinfo{author}{S.~Liu},
\newblock \bibinfo{title}{Fast cumulant method for probabilistic power flow considering the nonlinear relationship of wind power generation},
\newblock \bibinfo{journal}{IEEE Transactions on Power Systems} \bibinfo{volume}{35} (\bibinfo{year}{2020}) \bibinfo{pages}{2537--2548}.
\bibitem[{Hu and Wang(2006)}]{linear_pf2}
\bibinfo{author}{Z.~Hu}, \bibinfo{author}{X.~Wang},
\newblock \bibinfo{title}{A probabilistic load flow method considering branch outages},
\newblock \bibinfo{journal}{IEEE Transactions on Power Systems} \bibinfo{volume}{21} (\bibinfo{year}{2006}) \bibinfo{pages}{507--514}.
\bibitem[{Schieferdecker and Huber(2009)}]{GMMreduction1}
\bibinfo{author}{D.~Schieferdecker}, \bibinfo{author}{M.~F. Huber},
\newblock \bibinfo{title}{Gaussian mixture reduction via clustering},
\newblock in: \bibinfo{booktitle}{2009 12th International Conference on Information Fusion}, \bibinfo{year}{2009}, pp. \bibinfo{pages}{1536--1543}.
\bibitem[{Crouse et~al.(2011)Crouse, Willett, Pattipati, and Svensson}]{GMMreduction2}
\bibinfo{author}{D.~F. Crouse}, \bibinfo{author}{P.~Willett}, \bibinfo{author}{K.~Pattipati}, \bibinfo{author}{L.~Svensson},
\newblock \bibinfo{title}{A look at gaussian mixture reduction algorithms},
\newblock in: \bibinfo{booktitle}{14th International Conference on Information Fusion}, \bibinfo{year}{2011}, pp. \bibinfo{pages}{1--8}.
\bibitem[{Jabr(2006)}]{SOCP1}
\bibinfo{author}{R.~Jabr},
\newblock \bibinfo{title}{Radial distribution load flow using conic programming},
\newblock \bibinfo{journal}{IEEE Transactions on Power Systems} \bibinfo{volume}{21} (\bibinfo{year}{2006}) \bibinfo{pages}{1458--1459}.
\bibitem[{Kayacık and Kocuk(2021)}]{SOCP2}
\bibinfo{author}{S.~E. Kayacık}, \bibinfo{author}{B.~Kocuk},
\newblock \bibinfo{title}{An misocp-based solution approach to the reactive optimal power flow problem},
\newblock \bibinfo{journal}{IEEE Transactions on Power Systems} \bibinfo{volume}{36} (\bibinfo{year}{2021}) \bibinfo{pages}{529--532}.
\bibitem[{Farivar and Low(2013)}]{SOCP3}
\bibinfo{author}{M.~Farivar}, \bibinfo{author}{S.~H. Low},
\newblock \bibinfo{title}{Branch flow model: Relaxations and convexification—part i},
\newblock \bibinfo{journal}{IEEE Transactions on Power Systems} \bibinfo{volume}{28} (\bibinfo{year}{2013}) \bibinfo{pages}{2554--2564}.
\bibitem[{Yu et~al.(2018)Yu, Yin, and Zhu}]{STGCN}
\bibinfo{author}{B.~Yu}, \bibinfo{author}{H.~Yin}, \bibinfo{author}{Z.~Zhu},
\newblock \bibinfo{title}{Spatio-temporal graph convolutional networks: A deep learning framework for traffic forecasting},
\newblock in: \bibinfo{booktitle}{Proceedings of the Twenty-Seventh International Joint Conference on Artificial Intelligence}, \bibinfo{year}{2018}, pp. \bibinfo{pages}{3634--3640}.

\end{thebibliography}

\end{document}